\documentclass[11pt]{amsart}
\usepackage{amsmath, amsfonts, amssymb, amsthm, geometry, graphics, graphicx, multicol, multirow, enumerate, float, floatflt, fullpage, tikz, outlines, url, tabularx, caption}
\usetikzlibrary{decorations.text}
\usepackage{amscd,mathrsfs}
\usepackage{hyperref}
\usepackage{color}
\usepackage{graphicx}
\usepackage{stmaryrd}  %% double square brackets
\usepackage{tikz-cd}
\usetikzlibrary{shapes.multipart,shapes, matrix, positioning, arrows,cd,tqft}
  
\let\oldtocsection=\tocsection
\let\oldtocsubsection=\tocsubsection
\renewcommand{\tocsection}[2]{\hspace{0em}\oldtocsection{#1}{#2}}
\renewcommand{\tocsubsection}[2]{\hspace{1em}\oldtocsubsection{#1}{#2}}

\newtheorem{theorem}{Theorem}

\newtheorem{prop}[theorem]{Proposition}

\newcommand\Z{{\mathbb{Z}}}

\newcommand{\oplusop}[1]{{\mathop{\oplus}\limits_{#1}}}

\def\lra{\longrightarrow}
\def\Hom{\mathrm{Hom}}
\def\End{\mathrm{End}}
\def\Tr{\mathrm{Tr}}

\def\Cob{\mathrm{Cob}} 
\def\Kob{\mathrm{Kob}}

\def\Cobtwo{\Cob_2}   %  2D cobordisms
\def\Cobal{\Cob_{\alpha}}     % 2D cobordisms  mod alpha
\def\Cobalp{\Cob_{\alpha}'} 
\def\Kobal{\Kob_{\alpha}}   %  Karoubi envelope  
\def\PCob{\mathrm{PCob}}
\def\PCobal{\PCob_{\alpha}}

\def\KPobal{\mathrm{PKob}_{\alpha}}
\def\Kar{\mathrm{Kar}}

\def\Pa{\mathrm{Pa}}   % partitions

\def\lra{\longrightarrow}

\def\kk{\mathbf{k}}  %% base field  

\def\dmod{\mathrm{-mod}}   % modules  
\def\delcat{\mathrm{Rep}(S_t)}   % deligne category 
\def\delcatn{\mathrm{Rep}(S_n)}   % deligne category 
\newcommand{\rseries}[1]{R \llbracket #1 \rrbracket}    % power series 

\def\R{\mathbb R}

\def\Z{\mathbb Z}
\def\N{\mathbb N}

\def\SS{\mathbb{S}}    % circle 

\let\emptyset\varnothing

\title{Bilinear pairings on two-dimensional cobordisms and generalizations of the Deligne category}
\author{Mikhail Khovanov}
 \address{Department of Mathematics, Columbia University, New York, NY 10027, USA}
 \email{\href{mailto:khovanov@math.columbia.edu}{khovanov@math.columbia.edu}}
\author{Radmila Sazdanovic}
\address{Department of Mathematics, North Carolina State University, Raleigh, NC 27696-8205, USA}
 \email{\href{mailto:rsazdan@ncsu.edu}{rsazdan@ncsu.edu}}

\date{July 24, 2020}

\begin{document}

\tikzstyle{every node}=[font=\large]
\tikzstyle{every path}=[line width=1pt]

\begin{abstract}
The Deligne category of symmetric groups is the  additive Karoubi closure of the partition category. It is semisimple for generic values of  the parameter $t$ while producing categories of representations of the symmetric group when 
modded out by the ideal of negligible morphisms when $t$ is  a non-negative integer. The partition category 
may be interpreted, following Comes, via a particular linearization of the category of two-dimensional oriented cobordisms. The Deligne category and  its semisimple quotients admit similar interpretations. This viewpoint coupled to the universal construction of two-dimensional topological theories leads to multi-parameter monoidal generalizations of the partition  and the Deligne categories, one for each rational function in one variable.  
\end{abstract}

\maketitle
\tableofcontents

%%%%%%%%%%%%%%%%%%%%%
%
%  Introduction
%
%%%%%%%%%%%%%%%%%%%%%

\section{Introduction}

\smallskip

The Deligne category $\delcat$ interpolates between the categories of finite-dimensional representations of the symmetric groups $S_n$, viewed as tensor categories, turning integer $n$ into an element $t$ of the ground field~\cite{D},~\cite{CO},~\cite[Section 9.12.1]{EGNO}. 

The Deligne category has a diagrammatic description, via the partition category $\Pa_t$, as the Karoubi envelope of the additive closure of  $\Pa_t$.  When $t=n$ is a non-negative integer, the Deligne category $\delcatn $ has a non-trivial ideal of negligible morphisms, and the quotient  by this ideal  is naturally equivalent to the tensor  category of  finite-dimensional representations  of the symmetric group.

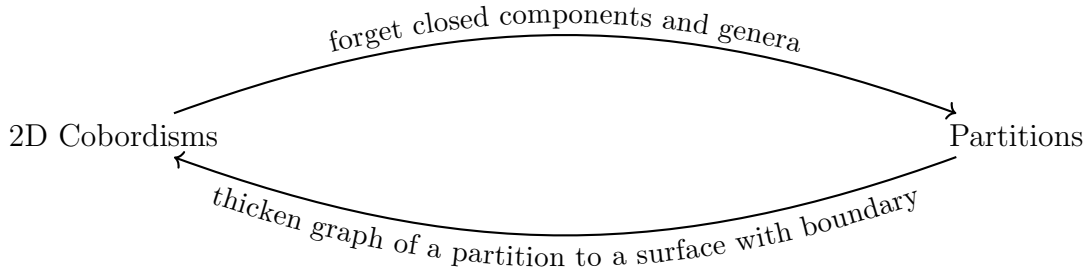
\begin{figure}[h]
\begin{center}
\begin{tikzpicture}
\path (0,0) node[centered](x) {2D Cobordisms} 
(12,0) node(y) {Partitions};

\def\myshift#1{\raisebox{-2.5ex}}

\draw [->,thick,postaction={decorate,decoration={text along path,reverse path,text align=center,text={|\myshift|thicken graph of a partition to a surface with boundary}}}] (y) to [bend right=-20]  (x);

\def\myshift#1{\raisebox{1ex}}
\draw [->,thick,postaction={decorate,decoration={text along path,text align=center,text={|\myshift|forget closed components and genera}}}]      (x) to [bend left=20] (y);

\end{tikzpicture}
\caption{Schematic correspondence between set partitions and 2D cobordisms.}
\label{fig_1_0}
\end{center}
\end{figure}

Diagrams commonly used to describe partitions~\cite{CO,C,HR,LS} can be  thickened to two-dimensional surfaces or cobordisms $S$ between unions  of circles. Circles appear as "thickenings" of points on which the partitions are formed. Vice versa, any 2D cobordism $S$ gives rise to  a partition  upon  ignoring  closed components of $S$ and the genus of each connected component with boundary. This informal correspondence  is depicted  in Figure~\ref{fig_1_0}.
Cobordisms boast higher variability than partitions, admitting components without boundary and allowing arbitrary genus of each component.

\begin{figure}[h]
\begin{center}
\begin{tikzpicture}

\begin{scope}
  \draw[style=thick] (-0.5,1.5)--(-0.5,2)--(0.5,1.5)--(0.5,-0.5)--(-0.5,0)--(-0.5,0.5);
  \draw[style=thick] (-0.8,1.4) to[out=0, looseness = 1.5,in=240] (-0,1.5); 
  \draw[style=thick] (-0.8,1.4) to[out=180, looseness = 4,in=180] (-0.8,0.6); 
    \draw[style=thick] (-0.8,0.6) to[out=0, looseness = 1.5,in=140] (-0,0.5);
  \draw[style=thick] (-1.3,1) to[out=90, looseness = 1.1,in=90] (-0.9,1); 
  \draw[style=thick] (-1.4,1.1) to[out=270, looseness = 0.9,in=270] (-0.8,1.1); 

\end{scope}

\begin{scope}[shift={(2,0)}]
\node at (-1,1) {=};
\draw[style=thick] (-0.5,2)--(0.5,1.5)--(0.5,-0.5)--(-0.5,0)--(-0.5,2);
            \end{scope}
         
 \begin{scope}[shift={(5,1)}]
         \draw[radius=0.7](0,0) circle;
         \draw[style=thick, dashed] (-0.7,0) to[out=90, looseness = 0.5,in=90] (0.7,0);
            \draw (-0.7,0) to[out=-90, looseness = 0.5,in=-90] (0.7,0);
         \node at (1.5,0) {\Large{$=t$}};
 \end{scope}
\end{tikzpicture}
\caption{Handle  removal and sphere evaluation skein relations on 2D cobordisms.}
\label{fig_2_0}
\end{center}
\end{figure}
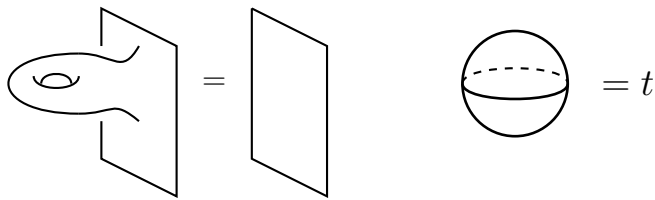

A precise connection between 2D cobordisms  and  the partition category was  pointed out by  Comes~\cite[Section 2.2]{C}:  modding out the cobordism category by the relations that adding a handle is the  identity and that a 2-sphere  evaluates to $t$, see  Figure~\ref{fig_2_0}, produces the partition category, with parameter $t$ corresponding to the 2-sphere. Comes used this observation to derive a set of defining  relations for the  partition category from that of  the cobordism category.

A family of  2-dimensional topological theories was recently introduced by one of the authors~\cite{Kh2}, based on Blanchet, Habegger, Masbaum and Vogel's \emph{universal construction}~\cite{BHMV}. It starts with an evaluation of closed oriented surfaces, which may be described by power series 
\begin{equation} \label{eq_Z_pow}
    Z_{\alpha}(T) = \alpha_0+\alpha_1 T + \alpha_2 T^2 +\dots =\sum_{n\ge 0} \alpha_n T^n \ \in \rseries{T},
\end{equation}
where  $R$ is a ground commutative ring or a field $\kk$, evaluating   a connected component of genus $g$ to $\alpha_g$. This  evaluation gives  rise to state spaces $A_{\alpha}(k)$ for collections of $k$ circles. It can then be extended to produce a category $\Cobal$ with objects non-negative integers  $n$ and  hom spaces $\Hom_{\Cobal}(n,m)$ being $R$-linear combinations of cobordisms between $n$ and $m$ circles modulo universal relations defined by the sequence $\alpha=(\alpha_0,\alpha_1,\dots)$, also see below.  

In the partition category $\Pa_t$, 
when two partitions are composed, each connected component of the composition 
that has no boundary points is evaluated to $t\in R$ and removed. In the  correspondence between 2D cobordisms  and partitions, these  components give rise to cobordisms of various genera,  which, in general, evaluate  to $\alpha_g$, where  $g$ is the  genus of the cobordism.  

To match the general $\alpha$-evaluation of 2D cobordisms to the evaluation by powers of $t$ in the partition and Deligne categories  specialize the  sequence $\alpha$ to 
\begin{equation}\label{seq_alpha_t}
    \alpha(t) = (t,t,t,\dots ), \  \   \alpha_g(t)=t \ \ \forall g \in  \Z_+.
\end{equation}
Then the additive Karoubi closure  of the category $\Cob_{\alpha(t)}$ is equivalent, as a tensor category, to the Deligne category 
\begin{equation}
    \Kar(\Cob_{\alpha(t)}^{\oplus}) \cong  \delcat  , 
\end{equation}
for $t\in \kk\setminus \Z_+$ (specializing to characteristic zero field $\kk$ as the ground ring). When  $t=n\in \Z_+$, the quotient of the Deligne category by the ideal $J_n$ of negligible morphisms  produces the category of finite-dimensional representations of the symmetric group  $S_n$, equivalent  to  the  above Karoubi  closure for $t=n$, 
\begin{equation}
    \Kar(\Cob_{\alpha(n)}^{\oplus}) \cong  \delcatn/J_n \cong \kk[S_n]\dmod. 
\end{equation}

This observation allows to generalize the Deligne category and its semisimple quotients by taking a more general sequence $\alpha$ of elements of $R$ and then forming tensor category $\Cobal$ and its additive Karoubi closure $\Kar(\Cobal^{\oplus})$, also denoted $\Kobal$. The latter is given  by  first allowing finite linear combinations of objects of $\Cobal$,  with suitably defined  hom spaces, and then adding all idempotents in endomorphism rings of these linear combinations as additional objects. 

When $R$ is a field $\kk$, it follows from \cite{Kh2} and goes back to a  theorem of Kronecker that hom spaces in $\Cobal$ are finite-dimensional iff the power series $Z_{\alpha}(T)$ in (\ref{eq_Z_pow}) can be represented as a rational function, 
\begin{equation}
    Z_{\alpha}(T) = \frac{P(T)}{Q(T)},
\end{equation}
where $P(T),Q(T)$ are coprime polynomials with coefficients  in  $\kk$. To  each such rational function we can assign an additive Karoubi-complete tensor (symmetric monoidal)  category 
\begin{equation} 
\Kobal:=\Kar(\Cobal^{\oplus})
\end{equation}
with finite dimensional hom spaces. This category is a natural generalization of the Deligne  category $\delcat$ for generic $t$ and of its semisimple quotients for $t=n\in \Z_+$. 
The Deligne category corresponds  to the rational function
\begin{equation}\label{eq_D_rat_fun}
    Z_{\alpha(t)} (T) = \frac{t}{1-T}  = t + tT + tT^2 + \dots  
\end{equation}
It should be  extremely interesting to extend various results and constructions related to the Deligne category and its semisimple quotients to this  large family of tensor categories $\Kobal$ (as well as  categories $\KPobal$ defined in Section~\ref{sec_gen_Deligne}) parametrized by rational functions. 

\vspace{0.1in} 

{\bf Acknowledgments:}  M.K. was partially supported by the NSF grant  DMS-1807425 while working on this paper. 
R.S. was partially supported by NSF grant DMS-1854705.

%%%%%%%%%%%%%%%%%%%%%%
% Linearization
%%%%%%%%%%%%%%%%%%%%%%

\section{Category of two-dimensional cobordisms and  its linearization categories}\label{sec_cat_lin}

\smallskip 

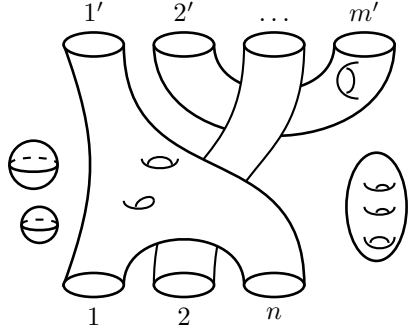
\begin{figure}[h]
\begin{center}
\begin{tikzpicture}
\begin{scope}[scale=0.8,shift={(0,0)}]

%left sphere bottom
\draw[radius=0.3](-0.9,1) circle;
         \draw[style=thick, dashed] (-1.2,1) to[out=90, looseness = 0.5,in=90] (-0.6,1);
            \draw (-1.2,1) to[out=-90, looseness = 0.5,in=-90] (-0.6,1);

% leftspehere top
\draw[radius=0.4](-1,2) circle;
         \draw[style=thick, dashed] (-1.4,2) to[out=90, looseness = 0.5,in=90] (-0.6,2);
            \draw (-1.4,2) to[out=-90, looseness = 0.5,in=-90] (-0.6,2);

%return
 \draw (1,4) to (2,4) to[out=270, looseness =1.3,in=270] (4,4) to (5,4) to [out=270, looseness =1.3,in=270]  (1,4);
 
%cylinder 
\draw[style=thick, fill=white]  (2.5,4) to (3.5,4) to [out=270, looseness =1.3,in=90] (2,0) to (1,0) to [out=90, looseness =1.3,in=270]  (2.5,4) ;

 %pants
 \draw[fill=white] (-0.5,4) to (0.5,4) to[out=270, looseness =1.3,in=90] (3.5,0) to (2.5,0) to [out=90, looseness =1.3,in=90]  (0.5,0)  to (-0.5,0) to[out=70, looseness =1.1,in=-70] (-0.5,4);
 
%genus
\draw[style=thick] (1.3,2) to[out=90, looseness = 1.1,in=90] (0.9,2); 
  \draw[style=thick] (1.4,2.1) to[out=270, looseness = 0.9,in=270] (0.8,2.1); 

\draw[style=thick] (1,1.4) to[out=90, looseness = 1.1,in=90] (0.7,1.3); 
  \draw[style=thick] (1,1.4) to[out=270, looseness = 0.9,in=270] (0.5,1.4);

\draw[style=thick] (4.2,3.6) to[out=0, looseness = 1.1,in=0] (4.2,3.2); 
  \draw[style=thick] (4.4,3.7) to[out=180, looseness = 2,in=180] (4.4,3.1);

%all genus to the right
 \draw[fill=white] (4.7, 1.3) ellipse (0.5cm and 0.9cm);
 
 \draw[style=thick] (4.7,1.6) to[out=90, looseness = 1.1,in=90] (4.9,1.6); 
  \draw[style=thick] (4.5,1.7) to[out=270, looseness = 0.9,in=270] (5,1.7);
  
   \draw[style=thick] (4.7,1.2) to[out=90, looseness = 1.1,in=90] (4.9,1.2); 
  \draw[style=thick] (4.5,1.3) to[out=270, looseness = 0.9,in=270] (5,1.3);
  
    \draw[style=thick] (4.5,0.8) to[out=270, looseness = 1.3,in=270] (5,0.8);
   \draw[style=thick] (4.6,0.7) to[out=90, looseness = 1.1,in=90] (4.9,0.7);
   
 \draw[fill=white] (0,4) ellipse (0.5cm and 0.2cm);
 \draw[fill=white](1.5,4) ellipse (0.5cm and 0.2cm);
 \draw[fill=white] (3,4)ellipse (0.5cm and 0.2cm);
 \draw[fill=white] (4.5,4)ellipse (0.5cm and 0.2cm);
  \node[above, font={\small}] at (0,4.2)  {$1'$};
  \node[above, font={\small}] at (1.5,4.2)  {$2'$};
    \node[above, font={\small}] at (3,4.2)  {$\ldots$};
\node[above, font={\small}] at (4.5,4.2)  {$m'$};
  
 \draw[fill=white] (0,0) ellipse (0.5cm and 0.2cm);
 \draw[fill=white](1.5,0) ellipse (0.5cm and 0.2cm);
 \draw[fill=white](3,0) ellipse (0.5cm and 0.2cm);
 \node[below, font={\small}] at (0,-0.2)  {$1$};
 \node[below, font={\small}] at (1.5,-0.2)  {$2$};
 \node[below, font={\small}] at (3,-0.2)  {$n$};
\end{scope}  
\end{tikzpicture}
\caption{A morphism in $\Cobtwo$. The  cobordism is not  embedded anywhere, so overlaps of components do not carry any  information and  can be reversed. We label top circles by  $1',2',\dots,m'$ and bottom circles by $1,2,\dots, n$. In this example $m=4$ and $n=3$.}
\label{fig_1_1}
\end{center}
\end{figure}

{\bf Category $\Cobtwo$.}
Consider the symmetric monoidal  category of 2-dimensional oriented cobordisms. We use the skeletal version of this category (one object in each isomorphism class), denoted  $\Cobtwo$. Its  objects are non-negative  integers  $n\in\Z_+  =\{0,1,2,\dots\}$ and morphisms from $n$ to $m$ are diffeomorphism classes rel boundary of compact oriented $2$-manifolds $S$ with a fixed diffeomorphism
\begin{equation}
    \partial S \cong (-\sqcup_n \SS^1)\sqcup (\sqcup_m \SS^1),
\end{equation}
where $\SS^1$  is   the oriented circle. In other words, the boundary of $S$ is separated into the \emph{bottom} and \emph{top} boundary, and identified, correspondingly, with disjoint unions of $n$ and $m$ circles. Composition is given by concatenation. Cobordisms may have connected components with  no boundary. A  example of a morphism (cobordism)  from $3$ to $4$  is given in  Figure~\ref{fig_1_1}. 

Morphisms from $n$  to $m$ in $\Cobtwo$ can be enumerated as follows. A morphism $x$ may have some number of closed components  of  various genera. Counting these components gives a sequence $cl(x)=(a_0,a_1,\dots, 0,0,\dots )$, where $a_k$ is the  number of closed components of genus $k$ in $x$. All  but finitely many terms  in the sequence $cl(x)$ are zero. Connected components with boundary provide a decomposition of the set of $n+m$ boundary circles into nonempty subsets, where circles from the same subset are the boundaries  of the same connected component. Furthermore, each such component has genus zero or higher, which counts the number of handles of the component.

Denote by $D^m_n$ the set  of decompositions of  $n+m$ circles. A morphism $x\in \Hom_{\Cobtwo}(n,m)$ can be described by 
a decomposition $\lambda\in D^m_n$, an assignment of a nonnegative integer (genus or  number of handles)  to each set in the decomposition $\lambda$ and a choice of a sequence $cl(x)$ as  above describing genera of closed components of $x$. 

Let us label bottom circles $1,\dots,n$ and  top circles $1',\dots, m'$, from left to right. Cobordism $x$ induces a decomposition of  the set 
\begin{equation}\label{eq_set_n_m}
\N_n^m := \{1,2,\dots, n,1',2',\dots, m'\} .  
\end{equation} 
 For the cobordism  $x$ in Figure~\ref{fig_1_1}  we  have $n=3$, $m=4$, the  set is $\N^4_3=\{1,2,3,1',2',3',4'\}$, and its subsets  corresponding  to components  with boundary are $\{1,3,1'\}$, $\{2,3'\}$, and $\{2',4'\}$.  These  components  have  genera   $2,0,1$, correspondingly.  The sequence  $cl(x)=(2,0,0,1,0,\dots)$, since $x$ has two closed components  of genus  zero (2-spheres) and  one  component of  genus three.  

\vspace{0.07in} 

{\bf Category $R\Cobtwo$.}
Fix a commutative ring  $R$  and consider  pre-additive $R$-linear category  $R\Cobtwo$ freely generated  by  $\Cobtwo$. 
It has the same  objects $n$ as $\Cobtwo$, and morphisms from $n$ to  $m$ in $R\Cobtwo$ are linear combinations of morphisms from $n$ to $m$ in $\Cobtwo$ with  coefficients in $R$ and with composition induced  from  that  in $\Cobtwo$. 
One can think of this construction, for an  arbitrary category $C$, as  analogous to passing from a   group  $G$ to its group algebra $R[G]$ or from a  semigroup $G$ to its semigroup algebra. 
It results in an idempotented ring $RC$ with a  collection of mutually  orthogonal idempotents (corresponding to identity morphisms), one for each object of $C$, as  a substitute for  the unit element, see~\cite{KS2}. 

\vspace{0.07in} 

{\bf Category $\Cobalp$ for a sequence $\alpha$.}
A more interesting  category is obtained if we  choose an infinite  sequence of elements 
\begin{equation}\label{eq_alpha}
\alpha=(\alpha_0,\alpha_1, \alpha_2, \dots )
\end{equation}
of $R$ and  evaluate each closed component of genus $k$ to $\alpha_k$. For a closed surface $S$ denote 
\begin{equation} 
\alpha(S) = \prod_{k\ge 0} \alpha_k^{c_k},
\end{equation} 
where $c_k$ is the number of components  of  $S$ of genus $k$. 
The resulting monoidal category, denoted  $\Cob'(\alpha)$ or $\Cobalp$, has the same  objects $n$ as  the  earlier  categories.  
A morphism from $n$ to $m$ in $\Cobalp$ is an  $R$-linear combination of cobordisms from $n$ circles  to  $m$ circles in $\Cobtwo$ without closed components. 
Composition is given  by concatenation followed by evaluating each closed component of genus $k$ to  $\alpha_k.$  We  can informally refer to $\Cobalp$ as the \emph{$\alpha$-prelinearization}  of the category $\Cobtwo$. 

There is the obvious "evaluation" or "reduction" functor $R\Cobtwo \lra \Cobalp$, which is the identity on objects, that evaluates (or reduces) each closed component of genus $k$ to  $\alpha_k$. The hom space $\Hom_{\Cobalp}(n,m)$ is a free $R$-module with a basis of cobordisms without closed components. Basis elements are parametrized by partitions in $D^m_n$ with a non-negative integer (genus) assigned to each part of the partition.

\vspace{0.07in} 

\emph{Remark:} Object $0$ associated to the empty $1$-manifold $\emptyset_1$ is the unit object of monoidal categories $\Cobtwo, R\Cobtwo$ and $\Cobalp$.  Commutative monoid of   endomorphisms  $\End_{\Cobtwo}(0)$ is freely generated  by isomorphism classes of closed oriented connected surfaces, one for each genus $g\ge 0$, and can  be  identified with the free abelian monoid on these generators.  Commutative rings 
\begin{equation*} 
\End_{R\Cobtwo}(0)\cong \End_{\Cobalp}(0)\cong R[  \End_{\Cobtwo}(0)]
\end{equation*}
are the  semigroup algebras  of  that monoid.  

\vspace{0.1in} 

{\bf Generating functions and the bilinear form.}
Sequence  $\alpha$ is conveniently encoded by the generating function 
\begin{equation}\label{eq_gen_fun} 
    Z_{\alpha}(T) = \sum_{n\ge 0}\alpha_n T^n \in \rseries{T}.
\end{equation}
For a closely related construction see~\cite{Kh2}, where to such $Z_{\alpha}(T)$ there is associated a family of $R$-modules  $A_{\alpha}(n)$, for each $n\ge 0$, constructed via a bilinear form on the space of linear combinations of oriented 2-manifolds with boundary the disjoint union of $n$ circles $\sqcup_n \SS^1$. Namely, one considers the free $R$-module $\mathrm{Fr}(n)$ with a basis $\{[S]\}_S$ of oriented compact surfaces $S$ with $\partial S \cong \sqcup_n  \SS^1$ (with the diffeomorphism fixed). On $\mathrm{Fr}(n)$  there is an $R$-bilinear form $(,)_n$ given on pairs of generators $S_1,S_2$ by gluing the two surfaces along the common boundary and evaluating via $\alpha$: 
\begin{equation}\label{eq_bilin_f}
    ([S_1],[S_2])_n= \alpha ((-S_1)\sqcup_{\partial} S_2)
\end{equation}
The state space of $n$ circles is the quotient of $\mathrm{Fr}(n)$ by  the kernel of this  bilinear form: 
\begin{equation}
    A_{\alpha}(n)\ := \ \mathrm{Fr}(n) / \mathrm{ker}((,)_n).
\end{equation}

This collection of $R$-modules is naturally a representation of the category $\Cobalp$, when the latter is viewed as an idempotented $R$-algebra  with a system of mutually-orthogonal idempotents $\{1_n\}_{n\in \N}$. 
Namely, to the object $n$ of $\Cobalp$ associate the $R$-module $A_{\alpha}(n)$. To a morphism given by a cobordism $x\in \Cobtwo$ from $n$ to $m$  associate an $R$-module map 
\begin{equation}\label{eq_x_alpha}
x_{\alpha}: A_{\alpha}(n)\lra A_{\alpha}(m),
\end{equation} 
obtained directly from the construction in~\cite{Kh2}, via the evaluations of $x$  capped off by various  oriented surfaces with $n$ and $m$ circles as the boundary. These morphisms over all $x\in  \Cobtwo$ provide a  representation of  $\Cobtwo$ and  $\Cobalp$ on the direct sum of $R$-modules
\begin{equation} \label{eq_A_rep}
    A_{\alpha} : = \oplusop{n\ge 0} A_{\alpha}(n).
\end{equation}
Monoidal structures on $\Cobtwo$ and  $\Cobalp$ are not used in this construction. $A_{\alpha}$ is a representation of  the idempotented $R$-algebra underlying category $\Cob'_{\alpha}$, in the sense of~\cite{KS1,KS2}. 

\vspace{0.1in} 

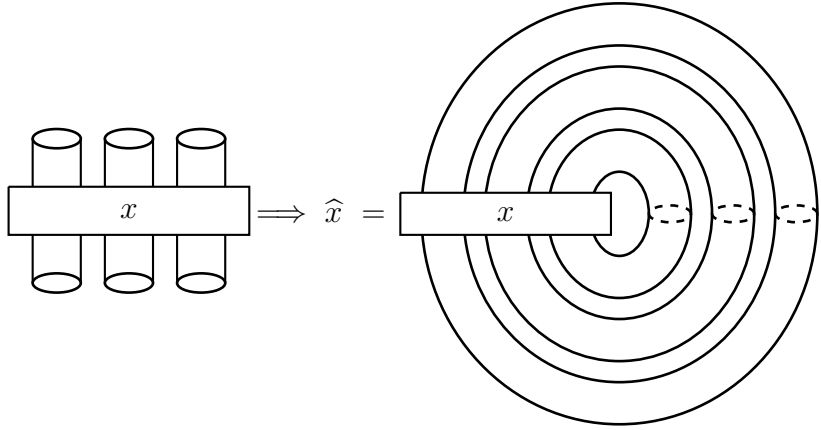
\begin{figure}[h]
\begin{center}
\begin{tikzpicture}[scale=0.8]
         \begin{scope}[scale=0.8,shift={(0,0)}]
\draw (1,2) ellipse (0.5cm and 0.2cm) ;         
 \draw (2.5,2) ellipse (0.5cm and 0.2cm);
 \draw (4,2) ellipse (0.5cm and 0.2cm); 
\draw (1,-1) ellipse (0.5cm and 0.2cm);         
 \draw (2.5,-1) ellipse (0.5cm and 0.2cm); 
 \draw (4,-1) ellipse ((0.5cm and 0.2cm);  

\draw[thick] (0.5,2)--(0.5, -1); 
\draw[thick] (1.5,2)--(1.5, -1);
\draw[thick] (2,2)--(2, -1);
\draw[thick] (3,2)--(3, -1);
\draw[thick] (3.5,2)--(3.5, -1);
\draw[thick] (4.5,2)--(4.5, -1);

\node at (2.5,0.5) {$x$};
\draw[style=thick, fill=white] (0,1)--(5,1)--(5,0)--(0,0)--(0,1);
\node at (2.5,0.5) {$x$};
\node at (6.4,0.5) {$\implies \widehat{x}\ = $ };
            \end{scope}
            
     \begin{scope}[scale=0.7,shift={(9.3,0)}]
 \draw (5.2,.5) ellipse (0.7cm and 1cm); 
  \draw (5.2,.5) ellipse (1.7cm and 2cm); 
  
  \draw (5.2,.5) ellipse (2.2cm and 2.5cm);
  \draw (5.2,.5) ellipse (3.2cm and 3.5cm);
  
  \draw (5.2,.5) ellipse (3.7cm and 4cm);
  \draw (5.2,.5) ellipse (4.7cm and 5cm);
  
\draw[style=thick, fill=white] (0,1)--(5,1)--(5,0)--(0,0)--(0,1);
\node at (2.5,0.5) {$x$};

 \draw[dashed] (9.4,0.5) ellipse (0.5cm and 0.2cm);
 \draw[dashed] (7.9,0.5) ellipse (0.5cm and 0.2cm); 
 \draw[dashed](6.4,0.5) ellipse (0.5cm and 0.2cm);  

            \end{scope}            
\end{tikzpicture}
\caption{Closing up a linear  combination $x$ of  $(n,n)$ cobordisms into a linear combination $\widehat{x}$ of closed  cobordisms.}
\label{fig_2_1}
\end{center}
\end{figure}

{\bf  Category $\Cobal$ as a quotient by negligible morphisms.}
The action of $\Cobalp$ can be quotiented down  to a smaller category. Categories $R\Cobtwo$ and $\Cobalp$ admit trace maps. Namely, given an element $x\in  \Hom(n,n)$, a finite linear combination of cobordisms with $n$ bottom and $n$ top circles, close up opposite circles $i$ and $i'$, $1\le i\le n$ by  annuli to  get  a linear combination of closed cobordisms $\widehat{x}$ and then evaluate the result via $\alpha$:
\begin{equation}
    \Tr(x) = \alpha(\widehat{x}) \  \in  R, 
\end{equation}
see Figure~\ref{fig_2_1}.  

For morphisms $x\in\Hom(n,m)$ and $y\in\Hom(m,n)$ we have  $\Tr(xy)=\Tr(yx)$. 

\vspace{0.07in} 
A morphism $x\in \Hom(n,m)$ in the category $\Cobalp$ (or in $R\Cobtwo$) is called  \emph{negligible} if for  any $y\in \Hom(m,n)$ the trace $\Tr(yx)=0$. Denote by $J(n,m)\subset \Hom(n,m)$ the subset of negligible morphisms from $n$ to $m$. 
This subset is an $R$-submodule of $\Hom(n,m)$, and the  union of 
$J(n,m)$, over all $n,m\ge 0$, is the tensor ideal $J_{\alpha}$ of $\Cobalp$. Define the category $\Cobal$ to be the quotient of $\Cobalp$ by this ideal,
\begin{equation}
    \Cobal \ :=  \ \Cobalp / J_{\alpha}.
\end{equation}
This category has  objects  $n\in\Z_+$, and 
\begin{equation}
    \Hom_{\Cobal}(n,m) \ = \ \Hom_{\Cobalp}(n,m)/J(n,m).
\end{equation}
Category $\Cobal$ is an $R$-linear tensor category with duals and a   non-degenerate trace: for  any $x\in\Hom_{\Cobal}(n,m)$, $x\not= 0$, there is  $y\in \Hom_{\Cobal}(n,m)$ such that  $\Tr(yx)\not= 0$. For  information about ideals of negligible morphisms and corresponding quotient categories we refer the reader to~\cite{EO,BW}. 

Starting with the category $R\Cobtwo$ instead of $\Cobalp$ in this  construction will result  in the quotient category isomorphic to  $\Cobal$. 

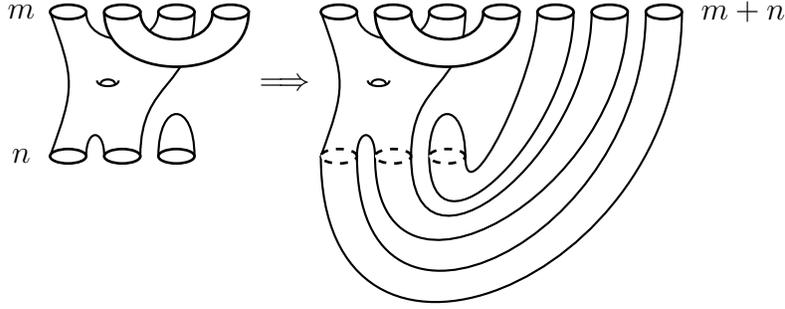
\begin{figure}[h]
\begin{center}

\begin{tikzpicture}[scale=0.6]

\draw [use as bounding box, white] (-2,4.5)-| (16,4.5)-| (16,-4)-| (-2,-4)-| (-2,4.5);
\begin{scope}[scale=0.8,shift={(0,0)}]
\node at (-1.3,4) {$m$};
\node at (-1.3,0) {$n$};
\node at (6,2) {$\implies$};
 \draw[style=thick] (-0.5,4) to (0.5,4) to[out=270, looseness =1.3,in=270] (2.5,4)to (3.5,4) to [out=270, looseness =1.3,in=90] (2,0) to (1,0) to [out=90, looseness =4,in=90] (0.5,0) to (-0.5,0) to [out=70, looseness =1.3,in=-70] (-0.5,4) ;

 \draw[style=thick] (1.3,2) to[out=90, looseness = 1.1,in=90] (0.9,2); 
  \draw[style=thick] (1.4,2.1) to[out=270, looseness = 0.9,in=270] (0.8,2.1); 
  
\draw[style=thick] (2.5,0)  to[out=90, looseness =4,in=90] (3.5,0);

 \draw[fill=white] (1,4) to (2,4) to[out=270, looseness =1.3,in=270] (4,4) to (5,4) to [out=270, looseness =1.3,in=270]  (1,4);

 \draw[fill=white] (0,4) ellipse (0.5cm and 0.2cm);
 \draw[fill=white](1.5,4) ellipse (0.5cm and 0.2cm);
 \draw[fill=white] (3,4)ellipse (0.5cm and 0.2cm);
 \draw[fill=white] (4.5,4)ellipse (0.5cm and 0.2cm);

 \draw[fill=white] (0,0) ellipse (0.5cm and 0.2cm);
 \draw[fill=white](1.5,0) ellipse (0.5cm and 0.2cm);
 \draw[fill=white](3,0) ellipse (0.5cm and 0.2cm);
\end{scope}   

\begin{scope}[scale=0.8,shift={(7.5,0)}]
 \draw[style=thick] (-0.5,4) to (0.5,4) to[out=270, looseness =1.3,in=270] (2.5,4)to (3.5,4) to [out=270, looseness =1.3,in=90] (2,0) to (1,0) to [out=90, looseness =4,in=90] (0.5,0) to (-0.5,0) to [out=70, looseness =1.3,in=-70] (-0.5,4) ;

 \draw[style=thick] (1.3,2) to[out=90, looseness = 1.1,in=90] (0.9,2); 
  \draw[style=thick] (1.4,2.1) to[out=270, looseness = 0.9,in=270] (0.8,2.1); 
  
\draw[style=thick] (2.5,0)  to[out=90, looseness =4,in=90] (3.5,0);

 \draw[fill=white] (1,4) to (2,4) to[out=270, looseness =1.3,in=270] (4,4) to (5,4) to [out=270, looseness =1.3,in=270]  (1,4);

 \draw[fill=white] (0,4) ellipse (0.5cm and 0.2cm);
 \draw[fill=white](1.5,4) ellipse (0.5cm and 0.2cm);
 \draw[fill=white] (3,4)ellipse (0.5cm and 0.2cm);
 \draw[fill=white] (4.5,4)ellipse (0.5cm and 0.2cm);
  \draw[fill=white](6,4) ellipse (0.5cm and 0.2cm);
 \draw[fill=white] (7.5,4)ellipse (0.5cm and 0.2cm);
 \draw[fill=white] (9,4)ellipse (0.5cm and 0.2cm);
 
 \draw[style=thick] (5.5,4)  to[out=270, looseness =1,in=270] (3.5,0);
  
 \draw[style=thick] (6.5,4)  to[out=270, looseness =1.5,in=270] (2.5,0);
 
  \draw[style=thick] (7,4)  to[out=270, looseness =1.6,in=270] (2,0);
  
  \draw[style=thick] (8,4)  to[out=270, looseness =1.6,in=270] (1,0);
 
 \draw[style=thick] (8.5,4)  to[out=270, looseness =1.8,in=270] (0.5,0);
  
 \draw[style=thick] (9.5,4)  to[out=270, looseness =1.8,in=270] (-0.5,0);
  
 \draw[fill=white, dashed] (0,0) ellipse (0.5cm and 0.2cm);
 \draw[fill=white, dashed](1.5,0) ellipse (0.5cm and 0.2cm);
 \draw[fill=white, dashed](3,0) ellipse (0.5cm and 0.2cm);
\node at (11.2,4) {$m+n$};
 \end{scope}     

\end{tikzpicture}
%\abovecaptionskip{-50pt}
\caption{Turning a morphism in $\Hom(n,m)$ into a morphism in $\Hom(0,n+m)$.}
\label{fig_2_2}
\end{center}
\end{figure}

The functor of modding out an $R$-linear tensor category with duals by the ideal of negligible morphisms is essentially the same operation as used in the universal construction~\cite{BHMV,Kh2}, where one mods out by  the kernel of the bilinear form. Thus, in the example above, there are  isomorphisms of $R$-modules  
\begin{equation}
\Hom_{\Cobal}(n,m)\cong \Hom_{\Cobal}(0,n+m)\cong  A_{\alpha}(n+m),
\end{equation}
with the first isomorphism given  by bending the $n$ bottom circles up, see Figure~\ref{fig_2_2}. 

 The space $J(0,n+m)$ of  negligible morphisms  in  $\Hom_{\Cobalp}(0,n+m)$ is exactly the kernel of the bilinear form  on $\Hom_{\Cobalp}(0,n+m)$ constructed via the formula (\ref{eq_bilin_f}) for $n+m$ boundary circles, implying the second isomorphism above.
 %\begin{equation}
%    \Hom_{\Cobal}(n,m) \cong A_{\alpha}(n+m).
%\end{equation}
It is easy to rewrite composition of morphisms in $\Cobal$ via these isomorphisms and suitable cobordism maps. 
   
\vspace{0.1in}    

We  refer to category $\Cobal$ as \emph{$\alpha$-linearization} of $\Cobtwo$  (and  of  related categories $R\Cobtwo$  and  $\Cobalp$). These categories are part of the package of the universal construction or pairing, see~\cite{BHMV,Kh2} and closely related~\cite{FKNSWW}, and can be defined in any dimension and in a variety of situations, given an evaluation of closed manifolds or similar objects (foams~\cite{Kh1,RW}, manifolds with embedded  submanifolds~\cite{FKNSWW,KR}, or other decorations). 
   
\vspace{0.1in}  

When commutative ring $R$ is a field $\kk$, it is observed in~\cite{Kh2} that the spaces $A_{\alpha}(n)$ are  finite-dimensional for all $n$ (equivalently, for  some $n\ge 1$)  iff the generating function (\ref{eq_gen_fun}) is a rational function in $T$.  Equvalently, representation (\ref{eq_A_rep}) of   
$\Cob_{\alpha}$ is locally finite-dimensional, in a similar sense, see~\cite{KS1,KS2}. The case when the function $Z_{\alpha}(T)$ is rational seems especially interesting, for many reasons. 

$R$-modules $A_{\alpha}(n)$  and maps between them induced by  cobordisms, see the discussion around (\ref{eq_x_alpha}),  define a representation of $\Cobal$  viewed  as an idempotented  ring 
\begin{equation} 
B_{\alpha} = \oplusop{n,m\ge 0} 1_m \Hom_{\Cobal}(n,m)1_n, 
\end{equation}
see~\cite{KS1,KS2}  for a  general discussion. On the corresponding representation
$A_{\alpha}$ in (\ref{eq_x_alpha}) 
idempotent $1_n$ acts as the projector 
onto $A_{\alpha}(n)$, and an element $x\in  \Hom_{\Cobal}(n,m) $ acts by  the corresponding  map  $A_{\alpha}(n)\lra A_{\alpha}(m)$. 
When $R$ is a field and $Z_{\alpha}(T)$ is rational, this representation is locally finite-dimensional. 

\vspace{0.1in} 

{\bf Additive closure and  the Karoubi envelope.} 
It is useful to consider the additive  Karoubi envelope $\Kob_{\alpha}$  of $\Cob_{\alpha}$. First 
form the finite additive closure  $\Cobal^{\oplus}$ of $\Cobal$ by taking formal finite direct sums of objects  $n$ of  $\Cobal$, and extending to morphisms in the obvious way. The additive closure has the zero object $\mathbf{0}$ different from the object $0$. The latter is  associated to the empty 1-manifold and comes from the corresponding object of $\Cobal$. Endomorphisms of the object $\mathbf{0}$ is the zero $R$-algebra, while $\End_{\Cobal^{\oplus}}(0)=\End_{\Cobal}(0)\cong R$. 

Denote the resulting category by $\Cobal^{\oplus}.$
Next, let $\Kobal$  be  the Karoubi  envelope of $\Cobal^{\oplus}.$ The six types of categories we've encountered so far are listed  below:
\begin{equation}\label{eq_many_cats} 
    \Cobtwo \lra R\Cobtwo\lra \Cobalp \lra \Cobal\lra \Cobal^{\oplus} \lra \Kobal.
\end{equation}
The first  arrow consists of allowing  $R$-linear combinations of cobordisms. In the second  arrow we  evaluate closed surfaces of genus $k$ to fixed elements $\alpha_k$ of  $R$, over all $k\ge 0$. We can refer to this procedure as  $\alpha$-prelinearization. The third arrow consists of modding  out  $\Cobalp$ by the ideal  of negligible morphisms. 

Like $\Cobalp$, category $R\Cobtwo$ also has the ideal $I_{\alpha}$ of negligible morphisms, via the trace given by $\alpha$. The  composition of the second and  third  arrows above can also  be  described as the quotient of $R\Cobtwo$ by this ideal 
\begin{equation} 
R\Cobtwo \lra R\Cobtwo/I_{\alpha} \cong \Cobal.
\end{equation} 

\vspace{0.1in} 

The fourth and the fifth arrows in (\ref{eq_many_cats}) are  fully faithful functors. The second, 
third and fourth categories are pre-additive, the fifth category is additive and the last category is additive and Karoubi-complete. All six categories are tensor (symmetric monoidal) and these five functors are monoidal. 

\vspace{0.1in}

%%%%%%%%%%%%%%%%%%%%%%%
%
% Partition category   
%
%%%%%%%%%%%%%%%%%%%%%%%

\section{Partition  category and the  Deligne category}

\smallskip 

Recall that $R$ is a commutative  ring. In the context of the partition category and the Deligne category ring $R$ is often taken to be a field $\kk$. Fix $t\in R$.   

\vspace{0.07in}

{\bf Partition category.}
Partition category $\Pa_t$ extends the notion of the partition algebra that originally appeared in Martin~\cite{M} and Jones~\cite{J}, see~\cite{HR,LS} for more information and references. 

Objects $n$ of the partition category are non-negative integers and morphisms from $n$  to $m$ are $R$-linear combinations of decompositions $D_n^m$ (also called partitions) of the set  
$\N_n^m$, see discussion around formula (\ref{eq_set_n_m}). 

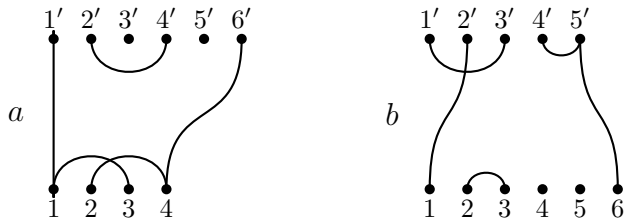
\begin{figure}[h]
\centering
\begin{tikzpicture}
 \begin{scope}[shift={(0,0)}]
\node at (-0.5,1) {$a$}; 
\draw[style=thick] (0,0) to[out=270,in=90] (0,2);
\draw[style=thick] (0,0) to[out=90, looseness = 1.5,in=90] (1,0);
\draw[style=thick] (0.5,0) to[out=90, looseness = 1.5,in=90] (1.5,0);
\draw[style=thick] (2.5,2) to[out=270, looseness = 1.5,in=90] (1.5,0);
\draw[style=thick] (0.5,2) to[out=270, looseness = 1.5,in=270] (1.5,2);

\draw[radius=.05, fill=black](0,0) circle;
\node[below, font={\small}] at (0,0)  {1};
\draw[radius=.05, fill=black](0.5,0) circle;
\node[below, font={\small}] at (0.5,0) {2};
\draw[radius=.05, fill=black](1,0) circle;
\node[below, font={\small}] at (1,0)  {3};
\draw[radius=.05, fill=black](1.5,0) circle;
\node[below, font={\small}] at (1.5,0) {4};
% \draw[radius=.1, fill=black](2,0) circle;
% \draw[radius=.1, fill=black](2.5,0) circle;

\draw[radius=.05, fill=black](0,2) circle;
\node[above, font={\small}] at (0,2) {$1'$};
\draw[radius=.05, fill=black](0.5,2) circle;
\node[above, font={\small}] at (0.5,2)  {$2'$};
\draw[radius=.05, fill=black](1,2) circle;
\node[above, font={\small}] at (1,2) {$3'$};
\draw[radius=.05, fill=black](1.5,2) circle;
\node[above, font={\small}] at (1.5,2) {$4'$};
\draw[radius=.05, fill=black](2,2) circle;
\node[above, font={\small}] at (2,2) {$5'$};
\draw[radius=.05, fill=black](2.5,2) circle;
\node[above, font={\small}] at (2.5,2)  {$6'$};
\end{scope}

 \begin{scope}[shift={(5,0)}]
\node at (-0.5,1) {$b$}; 
\draw[style=thick] (0,0) to[out=90, looseness = 1.5,in=270] (0.5,2);
\draw[style=thick] (0.5,0) to[out=90, looseness = 1.5,in=90] (1,0);
\draw[style=thick] (1.5,2) to[out=270, looseness = 1.5,in=270] (2,2);
\draw[style=thick] (2,2) to[out=270, looseness = 1.5,in=90] (2.5,0);
\draw[style=thick] (0,2) to[out=270, looseness = 1.5,in=270] (1,2);

\draw[radius=.05, fill=black](0,0) circle;
\node[below, font={\small}] at (0,0)  {1};
\draw[radius=.05, fill=black](0.5,0) circle;
\node[below, font={\small}] at (0.5,0) {2};
\draw[radius=.05, fill=black](1,0) circle;
\node[below, font={\small}] at (1,0)  {3};
\draw[radius=.05, fill=black](1.5,0) circle;
\node[below, font={\small}] at (1.5,0) {4};
\draw[radius=.05, fill=black](2,0) circle;
\node[below, font={\small}] at (2,0) {5};
\draw[radius=.05, fill=black](2.5,0) circle;
\node[below, font={\small}] at (2.5,0)  {$6$};

\draw[radius=.05, fill=black](0,2) circle;
\node[above, font={\small}] at (0,2) {$1'$};
\draw[radius=.05, fill=black](0.5,2) circle;
\node[above, font={\small}] at (0.5,2)  {$2'$};
\draw[radius=.05, fill=black](1,2) circle;
\node[above, font={\small}] at (1,2) {$3'$};
\draw[radius=.05, fill=black](1.5,2) circle;
\node[above, font={\small}] at (1.5,2) {$4'$};
\draw[radius=.05, fill=black](2,2) circle;
\node[above, font={\small}] at (2,2) {$5'$};

\end{scope}
\end{tikzpicture}

\caption{Partitions $a=\{\{1,3,1'\},\{2,4,6'\},\{2',4'\},\{3'\},\{5'\}\}\in P^6_4$ and $b=\{\{1,2'\},\{2,3\},\{4\},\{5\},\{6,4',5'\},\{1',3'\}\}\in P^5_6$. Notice multiple ways to display the same partition. For  the subset $\{1,3,1'\}$ we depicted edges $(1,3)$ and $(1,1')$. Another  possibility is to depict edges $(1,3)$  and $(3,1')$ or edges $(1,1')$ and $(3,1')$. In choosing a diagram for a partition It is natural to at least minimize the number of bottom-top edges, showing only one such edge for each subset that contains both bottom and  top points.}
\label{fig_1_2}

\end{figure}

Diagrammatically, partitions are often denoted by marking $n$ points on a horizontal line in the plane and $m$ points on a parallel line above it.  One connects these $n+m$ points by arcs, and connected components of the resulting graph are the parts of the partition. Intersections of arcs are  ignored.  A partition usually has more than one such diagram. For instance, if $\{1,3,1'\}$ is a part of  the partition, it can be described by two arcs $(1,3),(1,1')$ or two  arcs $(1,3),(3,1')$, or all three arcs, also see Figure~\ref{fig_1_5} left below demonstrating this indeterminacy.  
This diagrammatic description of partitions is standard in papers on the partition algebra and category, see for instance~\cite{LS}. 
Two examples of diagrammatic presentation are given in Figure~\ref{fig_1_2}. 

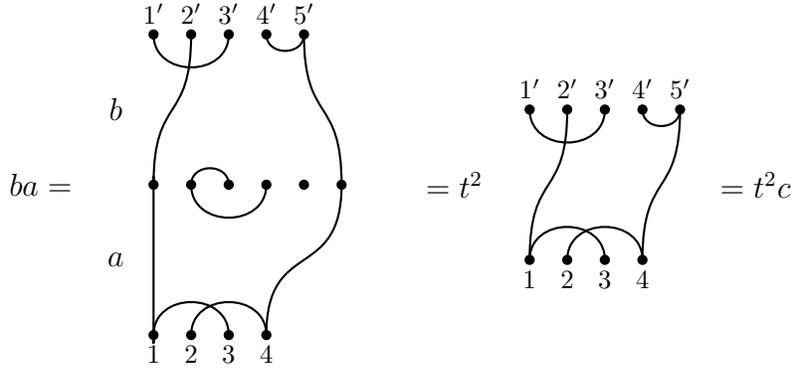
\begin{figure}[h]
\begin{center}
\begin{tikzpicture}
\node at (-1.5,2) {$ba=$}; 
          \begin{scope}[shift={(0,0)}]
\node at (-0.5,1) {$a$}; 
\draw[style=thick] (0,0) to[out=270,in=90] (0,2);
\draw[style=thick] (0,0) to[out=90, looseness = 1.5,in=90] (1,0);
\draw[style=thick] (0.5,0) to[out=90, looseness = 1.5,in=90] (1.5,0);
\draw[style=thick] (2.5,2) to[out=270, looseness = 1.5,in=90] (1.5,0);
\draw[style=thick] (0.5,2) to[out=270, looseness = 1.5,in=270] (1.5,2);

\draw[radius=.05, fill=black](0,0) circle;
\node[below, font={\small}] at (0,0)  {1};
\draw[radius=.05, fill=black](0.5,0) circle;
\node[below, font={\small}] at (0.5,0) {2};
\draw[radius=.05, fill=black](1,0) circle;
\node[below, font={\small}] at (1,0)  {3};
\draw[radius=.05, fill=black](1.5,0) circle;
\node[below, font={\small}] at (1.5,0) {4};
% \draw[radius=.1, fill=black](2,0) circle;
% \draw[radius=.1, fill=black](2.5,0) circle;

\draw[radius=.05, fill=black](0,2) circle;
% \node[above, font={\small}] at (0,2) {$1'$};
\draw[radius=.05, fill=black](0.5,2) circle;
% \node[above, font={\small}] at (0.5,2)  {$2'$};
\draw[radius=.05, fill=black](1,2) circle;
% \node[above, font={\small}] at (1,2) {$3'$};
\draw[radius=.05, fill=black](1.5,2) circle;
% \node[above, font={\small}] at (1.5,2) {$4'$};
\draw[radius=.05, fill=black](2,2) circle;
% \node[above, font={\small}] at (2,2) {$5'$};
\draw[radius=.05, fill=black](2.5,2) circle;
% \node[above, font={\small}] at (2.5,2)  {$6'$};
\end{scope}

 \begin{scope}[shift={(0,2)}]
\node at (-0.5,1) {$b$}; 
\draw[style=thick] (0,0) to[out=90, looseness = 1.5,in=270] (0.5,2);
\draw[style=thick] (0.5,0) to[out=90, looseness = 1.5,in=90] (1,0);
\draw[style=thick] (1.5,2) to[out=270, looseness = 1.5,in=270] (2,2);
\draw[style=thick] (2,2) to[out=270, looseness = 1.5,in=90] (2.5,0);
\draw[style=thick] (0,2) to[out=270, looseness = 1.5,in=270] (1,2);

\draw[radius=.05, fill=black](0,2) circle;
\node[above, font={\small}] at (0,2) {$1'$};
\draw[radius=.05, fill=black](0.5,2) circle;
\node[above, font={\small}] at (0.5,2)  {$2'$};
\draw[radius=.05, fill=black](1,2) circle;
\node[above, font={\small}] at (1,2) {$3'$};
\draw[radius=.05, fill=black](1.5,2) circle;
\node[above, font={\small}] at (1.5,2) {$4'$};
\draw[radius=.05, fill=black](2,2) circle;
\node[above, font={\small}] at (2,2) {$5'$};

\end{scope}

 \begin{scope}[shift={(5,1)}]
\node at (-1,1) {$=t^2$}; 
\draw[style=thick] (0,0) to[out=90, looseness = 1.5,in=270] (0.5,2);
\draw[style=thick] (0.5,0) to[out=90, looseness = 1.5,in=90] (1.5,0);
\draw[style=thick] (0,0) to[out=90, looseness = 1.5,in=90] (1,0);
\draw[style=thick] (1.5,2) to[out=270, looseness = 1.5,in=270] (2,2);
\draw[style=thick] (2,2) to[out=270, looseness = 1.5,in=90] (1.5,0);
\draw[style=thick] (0,2) to[out=270, looseness = 1.5,in=270] (1,2);

\draw[radius=.05, fill=black](0,0) circle;
\node[below, font={\small}] at (0,0)  {1};
\draw[radius=.05, fill=black](0.5,0) circle;
\node[below, font={\small}] at (0.5,0) {2};
\draw[radius=.05, fill=black](1,0) circle;
\node[below, font={\small}] at (1,0)  {3};
\draw[radius=.05, fill=black](1.5,0) circle;
\node[below, font={\small}] at (1.5,0) {4};

\draw[radius=.05, fill=black](0,2) circle;
\node[above, font={\small}] at (0,2) {$1'$};
\draw[radius=.05, fill=black](0.5,2) circle;
\node[above, font={\small}] at (0.5,2)  {$2'$};
\draw[radius=.05, fill=black](1,2) circle;
\node[above, font={\small}] at (1,2) {$3'$};
\draw[radius=.05, fill=black](1.5,2) circle;
\node[above, font={\small}] at (1.5,2) {$4'$};
\draw[radius=.05, fill=black](2,2) circle;
\node[above, font={\small}] at (2,2) {$5'$};
\node at (3,1) {$=t^2  c$}; 
\end{scope}
\end{tikzpicture}
\caption{Composition $ba = t^2  c$, where partition $c$ is shown  on the right,  $c=\{\{1,3,2'\},\{2,4,4',5'\},\{1',3'\} \}$. Coefficient $t^2$ comes from removing two connected components in the middle of $ba$ diagram that connect to neither bottom nor top points.}
\label{fig_1_3}
\end{center}
\end{figure}

Composition is given by concatenating diagrams, see Figure~\ref{fig_1_3}, and treating points in the middle that connect to bottom or top as 'pass through' points that vanish from the concatenation but are used before that to create the new partition. If there exist a connected component that consists entirely of points in the middle part of the diagram, it is removed and what's left is multiplied by $t$. This procedure is iterated until no such components are left. 

Composition is  then  extended bilinearly to $R$-linear combinations of partitions. The resulting $R$-linear category $\Pa_t$ is symmetric monoidal, with the  tensor product given on partitions by  placing their diagrams in parallel. 

\vspace{0.1in}

\begin{figure}[h]
\begin{center}
\begin{tikzpicture}
        
%a1
\begin{scope}[shift={(3,0)}]
% \node at (-0.5, 1.8) {(a)};
% \node at (2.5, 1) {$\rightarrow$};
\draw (2,0) to (1,1) to (0,0);
\draw (1,1) to (1,2);
\draw[radius=.05, fill=black](0,0) circle;
\draw[radius=.05, fill=black](2,0) circle;
\draw[radius=.05, fill=black](1,2) circle;
\end{scope}
%a2
\begin{scope}
\node at (-0.5, 1.8) {(a)};
\node at (2.5, 1) {$\rightarrow$};
\draw (2,0) to [out=155, looseness = 1.5,in=25](0,0); %bottom
\draw (0,0) to [out=40, looseness = 1.5,in=260](1,2);%left
\draw[radius=.05, fill=black](0,0) circle;
\draw[radius=.05, fill=black](2,0) circle;
\draw[radius=.05, fill=black](1,2) circle;
\draw[radius=.05, fill=black](0,0) circle;
\draw[radius=.05, fill=black](2,0) circle;
% \draw[radius=.05, fill=black](1,1) circle;
% \draw[radius=.05, fill=black](0,2) circle;
% \draw[radius=.05, fill=black](1,2) circle;
\draw[radius=.05, fill=black](1,2) circle;
\end{scope}

%b1
\begin{scope}[shift={(7,0)}]
\node at (-0.5, 1.8) {(b)};
 \node at (2.5, 1) {$\rightarrow$};
\draw (1,0) to [out=105, looseness = 1.5,in=310] (0,2); %right
\draw (1,0) to [out=75, looseness = 1.5,in=230](2,2);

\draw[radius=.05, fill=black](1,0) circle;
\draw[radius=.05, fill=black](2,2) circle;
% \draw[radius=.05, fill=black](1,1) circle;
% \draw[radius=.05, fill=black](0,2) circle;
% \draw[radius=.05, fill=black](1,2) circle;
\draw[radius=.05, fill=black](0,2) circle;
\end{scope}
%b2
\begin{scope}[shift={(10,0)}]
\draw (2,2) to (1,1) to (0,2);
\draw (1,1) to (1,0);
\draw[radius=.05, fill=black](1,0) circle;
\draw[radius=.05, fill=black](2,2) circle;
% \draw[radius=.05, fill=black](1,1) circle;
% \draw[radius=.05, fill=black](0,2) circle;
% \draw[radius=.05, fill=black](1,2) circle;
\draw[radius=.05, fill=black](0,2) circle;
\end{scope}
%----------------------------------
%c1
\begin{scope}[shift={(0,-3)}]
\node at (-0.5, 1.8) {(c)};
 \node at (2.5, 1) {$\rightarrow$};
\draw[dashed] (0,0) to (2,0);
\draw[dashed] (0,2) to (2,2);
\draw[radius=.05, fill=black](1,2) circle;
\end{scope}
%c2
\begin{scope}[shift={(3,-3)}]
\draw[dashed] (0,0) to (2,0);
\draw[dashed] (0,2) to (2,2);
\draw (1,2) to (1,1);
\node at (1,1) {$\ast$};
\draw[radius=.05, fill=black](1,2) circle;
\end{scope}

%d1
\begin{scope}[shift={(9,-1)}, rotate=180]
\node at (2.5, 0.2) {(d)};
 \node at (-0.5, 1) {$\rightarrow$};
\draw[dashed] (0,0) to (2,0);
\draw[dashed] (0,2) to (2,2);
\draw[radius=.05, fill=black](1,2) circle;
\end{scope}
%d2
\begin{scope}[shift={(12,-1)},rotate=180]
\draw[dashed] (0,0) to (2,0);
\draw[dashed] (0,2) to (2,2);
\draw (1,2) to (1,1);
\node at (1,1) {$\ast$};
\draw[radius=.05, fill=black](1,2) circle;
\end{scope}
%-----------------------------------
%e1
\begin{scope}[shift={(3,-6)}]
% \node at (-0.5, 1.8) {(a)};
% \node at (2.5, 1) {$\rightarrow$};
\draw (0.2,0) to (1.8,2);
\draw (0.2,2) to (1.8,0);
\draw[radius=.05, fill=black](0.2,0) circle;
\draw[radius=.05, fill=black](1.8,0) circle;
% \draw[radius=.05, fill=black](1,1) circle;
\draw[radius=.05, fill=black](0.2,2) circle;
% \draw[radius=.05, fill=black](1,2) circle;
\draw[radius=.05, fill=black](1.8,2) circle;
\draw[dashed] (0,0) to (2,0);
\draw[dashed] (0,2) to (2,2);
\end{scope}
%e2
\begin{scope}[shift={(0,-6)}]
\node at (-0.5, 1.8) {(e)};
\node at (2.5, 1) {$\rightarrow$};
\draw (0.2,0) to (1.8,2);
\draw (0.2,2) to (1.8,0);
\draw[radius=.05, fill=black](0.2,0) circle;
\draw[radius=.05, fill=black](1.8,0) circle;
% \draw[radius=.05, fill=black](1,1) circle;
\draw[radius=.05, fill=black](0.2,2) circle;
% \draw[radius=.05, fill=black](1,2) circle;
\draw[radius=.05, fill=black](1.8,2) circle;
\draw[dashed] (0,0) to (2,0);
\draw[dashed] (0,2) to (2,2);
\end{scope}

%f1
\begin{scope}[shift={(7,-6)}]
\node at (-0.5, 1.8) {(f)};
 \node at (2.5, 1) {$\rightarrow$};
\draw[dashed] (0,0) to (2,0);
\draw[dashed] (0,2) to (2,2);
\draw[radius=.05, fill=black](1,1) circle;
\end{scope}
%d2
\begin{scope}[shift={(10,-6)}]
\draw[dashed] (0,0) to (2,0);
\draw[dashed] (0,2) to (2,2);
\node at (1,1.5) {$\ast$};
\draw (1,1.5) to (1,0.5);
\node at (1,0.5) {$\ast$};
\end{scope}

\end{tikzpicture}
\caption{Conventional and Snyder's generators for  the partition category. Object 0 is shown by a dashed line without dots on it.  In (c) and (d) we indicated the loose end of a strand by  the $\ast$ symbol; other ways to depict the end are fine too. Likeng-Savage~\cite{LS} use a similar notation for the generators (c),(d). Element shown  in (f) is a suitable composition of generators (c) and (d) and evaluates to $t$ in the partition category. Top and  bottom dashed lines in the depiction of a diagram are optional and are not shown in the top row diagrams.}
\label{fig_1_4}
\end{center}
\end{figure}
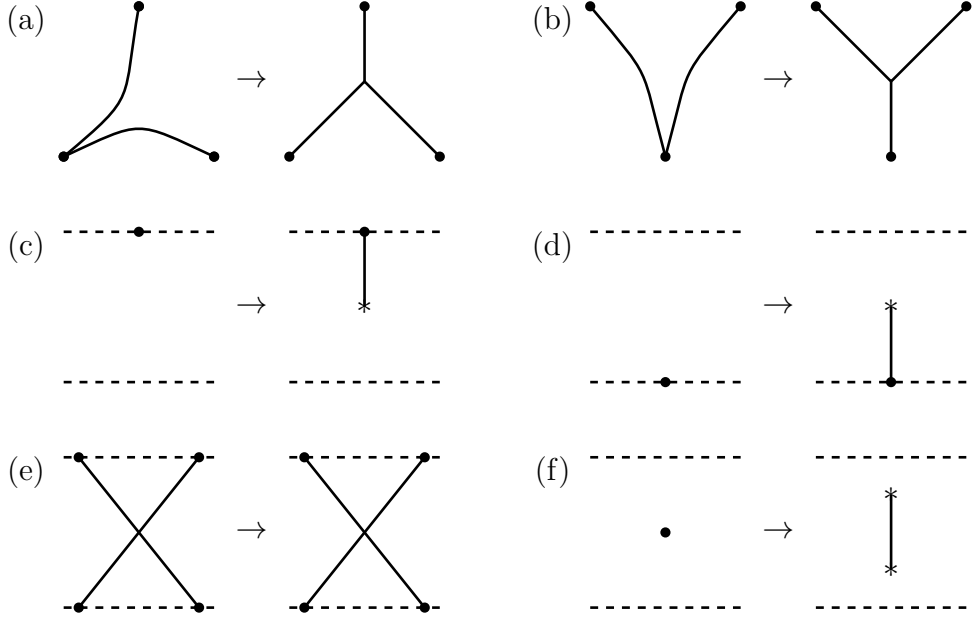

\begin{figure}[h]
\begin{center}
\begin{tikzpicture}[scale=0.9]
\begin{scope}
\node at (2.5, 1) {$=$};
\draw (1,2) to [out=-75, looseness = 1.5,in=130] (2,0); %right
\draw (0,0) to [out=50, looseness = 1.5,in=250](1,2);%left
\draw[radius=.05, fill=black](0,0) circle;
\draw[radius=.05, fill=black](2,0) circle;

\draw[radius=.05, fill=black](1,2) circle;
\end{scope}

\begin{scope}[shift={(3,0)}]
 \node at (2.5, 1) {$=$};
\draw (2,0) to [out=155, looseness = 1.5,in=25](0,0); %bottom
\draw (0,0) to [out=50, looseness = 1.5,in=250](1,2);%left
\draw[radius=.05, fill=black](0,0) circle;
\draw[radius=.05, fill=black](2,0) circle;
\draw[radius=.05, fill=black](1,2) circle;
\end{scope} 

\begin{scope}[shift={(6,0)}] \node at (2.5, 1) {$=$};
\draw (1,2) to [out=-75, looseness = 1.5,in=130] (2,0); %right
\draw (2,0) to [out=155, looseness = 1.5,in=25](0,0); %bottom
% \draw (2,0) to [out=50, looseness = 1.5,in=250](1,2);%left
\draw[radius=.05, fill=black](0,0) circle;
\draw[radius=.05, fill=black](2,0) circle;
\draw[radius=.05, fill=black](1,2) circle;
\end{scope}

\begin{scope}[shift={(9,0)}]
\draw (1,2) to [out=-75, looseness = 1.5,in=130] (2,0); %right
\draw (2,0) to [out=155, looseness = 1.5,in=25](0,0); %bottom
\draw (0,0) to [out=50, looseness = 1.5,in=250](1,2);%left
\draw[radius=.05, fill=black](0,0) circle;
\draw[radius=.05, fill=black](2,0) circle;
\draw[radius=.05, fill=black](1,2) circle;
\end{scope} 

\begin{scope}[shift={(13,0)}]
\draw (2,2) to  (1,1) to (0,2);
\draw (1,1) to (1,0);
\draw[radius=.05, fill=black](1,0) circle;
\draw[radius=.05, fill=black](2,2) circle;
% \draw[radius=.05, fill=black](1,1) circle;
% \draw[radius=.05, fill=black](0,2) circle;
% \draw[radius=.05, fill=black](1,2) circle;
\draw[radius=.05, fill=black](0,2) circle;
\end{scope}
\end{tikzpicture}
\caption{Multiple ways to depict generating morphism (a) in Figure~\ref{fig_1_4} versus unique up to isotopy diagram in the Snyder graphical calculus. }
\label{fig_1_5}
\end{center}
\end{figure}
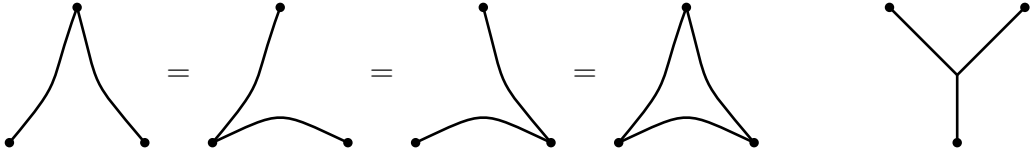

{\bf Noah Snyder's diagrammatics for the partition category.} Long time ago Noah Snyder~\cite{S} pointed out to one of us an alternative diagrammatics for the partition category, which will be treated in more detail in~\cite{H}. Figure~\ref{fig_1_4} shows conventional diagrammatics versus the Snyder diagrammatics for the standard generating morphisms of the partition category. One difference is the use of trivalent vertex to depict the morphism from $2$ to $1$ corresponding to the  partition $\{1,2,1'\}$ and the dual morphism from $1$ to $2$. This trivalent vertex as well as other configurations can be  freely rotated in the  plane. As in the usual diagrammatics,  one allows intersections of distinct parts of the partition, thinking of them as virtual intersections. 

Figure~\ref{fig_1_5} displays one benefit of the Snyder calculus: generating morphism (a) has essentially unique minimal presentation.

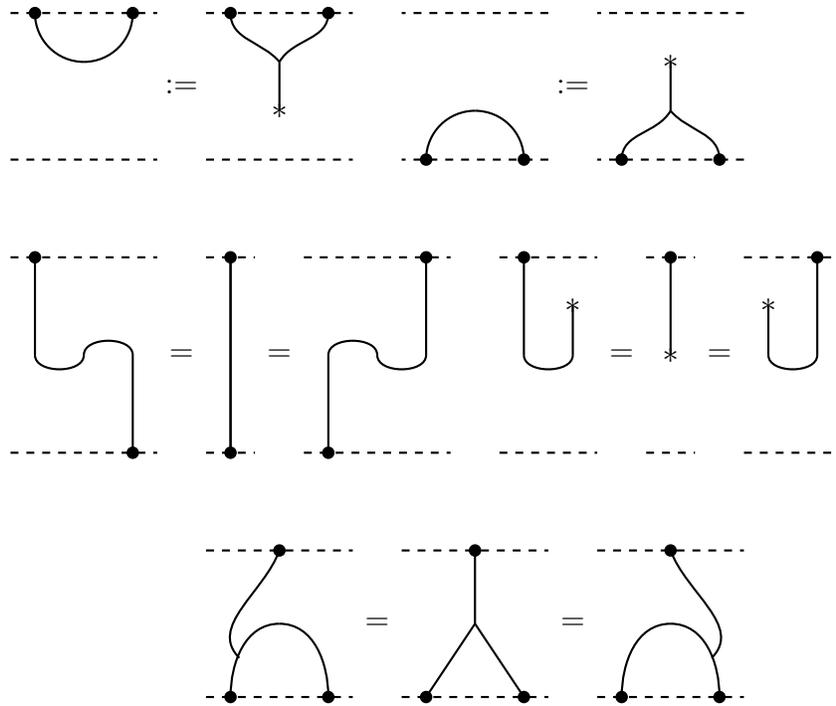
\begin{figure}[h]
\centering
\begin{tikzpicture} [baseline=0cm, scale=1.3]
%---------10-A1

\draw[style=thick, dashed](-0.25,1)--(1.25,1);
\draw[style=thick, dashed](-0.25,2.5)--(1.25,2.5);
\draw[style=thick] (0,2.5)to[out=270,in=180] (0.5,2);
\draw[style=thick] (0.5,2)to[out=0,in=-90] (1,2.5);
\draw[radius=.05, fill=black](0,2.5)circle;
\draw[radius=.05, fill=black](1,2.5)circle;
\node at (1.5,1.75) {:=};
%---------10-A2
\begin{scope}[shift={(2,0)}]
\draw[style=thick, dashed](-0.25,1)--(1.25,1);
\draw[style=thick, dashed](-0.25,2.5)--(1.25,2.5);
\draw[style=thick](0.5,2)--(0.5,1.5);
\node at (0.5,1.5) {$\ast$};
\draw[style=thick] (0,2.5)to[out=270,in=130] (0.5,2);
\draw[style=thick] (0.5,2)to[out=60,in=-90] (1,2.5);
\draw[radius=.05, fill=black](0,2.5)circle;
\draw[radius=.05, fill=black](1,2.5)circle;
\end{scope}
%---------10-B1
\begin{scope}[shift={(5,3.5)}, rotate=180]
\draw[style=thick, dashed](-0.25,1)--(1.25,1);
\draw[style=thick, dashed](-0.25,2.5)--(1.25,2.5);
\draw[style=thick] (0,2.5)to[out=270,in=180] (0.5,2);
\draw[style=thick] (0.5,2)to[out=0,in=-90] (1,2.5);
\draw[radius=.05, fill=black](0,2.5)circle;
\draw[radius=.05, fill=black](1,2.5)circle;
% \node at (1.75,1.5) {=:};
\end{scope}
%---------10-B2
\begin{scope}[shift={(7,3.5)},rotate=180]
\draw[style=thick, dashed](-0.25,1)--(1.25,1);
\draw[style=thick, dashed](-0.25,2.5)--(1.25,2.5);
\draw[style=thick](0.5,2)--(0.5,1.5);
\draw[style=thick] (0,2.5)to[out=270,in=130] (0.5,2);
\draw[style=thick] (0.5,2)to[out=60,in=-90] (1,2.5);
\node at (0.5,1.5) {$\ast$};
\draw[radius=.05, fill=black](0,2.5)circle;
\draw[radius=.05, fill=black](1,2.5)circle;
 \node at (1.5,1.75) {:=};
\end{scope}

%---------10-C1
\begin{scope}[shift={(0,-3)}]
\draw[style=thick, dashed](-0.25,1)--(1.25,1);
\draw[style=thick, dashed](-0.25,3)--(1.25,3);

\draw[style=thick](1,1)--(1,2);
\draw[style=thick](0,3)--(0,2);

\draw[style=thick] (1,2)to[out=90,in=90] (0.5,2);
\draw[style=thick] (0,2)to[out=270,in=270] (0.5,2);

\draw[style=thick](2,1)--(2,3);
\draw[style=thick, dashed](1.75,3)--(2.25,3);
\draw[style=thick, dashed](1.75,1)--(2.25,1);

\draw[radius=.05, fill=black](1,1)circle;
\draw[radius=.05, fill=black](0,3)circle;
\node at (1.5,2) {=};
\draw[style=thick](2,1)--(2,3);
\draw[radius=.05, fill=black](2,1)circle;
\draw[radius=.05, fill=black](2,3)circle;
\node at (2.5,2) {=};
\end{scope}
%---------10-C2
\begin{scope}[shift={(3,-3)}]
\draw[style=thick, dashed](-0.25,1)--(1.25,1);
\draw[style=thick, dashed](-0.25,3)--(1.25,3);
\draw[style=thick](1,3)--(1,2);
\draw[style=thick](0,1)--(0,2);
\draw[style=thick] (1,2) to [out=270,in=270](0.5,2);
\draw[style=thick] (0,2) to [out=90,in=90] (0.5,2);
\draw[radius=.05, fill=black](0,1)circle;
\draw[radius=.05, fill=black](1,3)circle;
\end{scope}

%---------10-D1
\begin{scope}[shift={(5,-3)}]
\draw[style=thick, dashed](-0.25,1)--(0.75,1);
\draw[style=thick, dashed](-0.25,3)--(0.75,3);

\draw[style=thick](0,3)--(0,2);
\draw[style=thick] (0,2)to[out=270,in=270] (0.5,2);
\draw[style=thick](0.5,2)--(0.5,2.5);
\draw[radius=.05, fill=black](0,3)circle;
\node at (0.5,2.5) {$\ast$};
\draw[radius=.05, fill=black](1.5,3)circle;
\node at (1,2) {=};
\draw[style=thick](1.5,2)--(1.5,3);
\node at (1.5,2) {$\ast$};
\draw[style=thick, dashed](1.25,3)--(1.75,3);
\draw[style=thick, dashed](1.25,1)--(1.75,1);
\node at (2,2) {=};
\end{scope}
%---------10-D2
\begin{scope}[shift={(7.5,-3)}]
\draw[style=thick, dashed](-0.25,1)--(0.75,1);
\draw[style=thick, dashed](-0.25,3)--(0.75,3);
\draw[style=thick](0.5,3)--(0.5,2);
\draw[style=thick](0,2.5)--(0,2);
\draw[style=thick] (0,2) to [out=270,in=270] (0.5,2);
\node at (0,2.5) {$\ast$};
\draw[radius=.05, fill=black](0.5,3)circle;
\end{scope}

%---------10-F1
\begin{scope}[shift={(2,-5.5)}]
\draw[style=thick, dashed](-0.25,1)--(1.25,1);
\draw[style=thick, dashed](-0.25,2.5)--(1.25,2.5);

\draw[style=thick] (0.5,1.75) to[out=180,in=90] (0,1);
\draw[style=thick] (0.5,1.75) to[out=0,in=90] (1,1);
\draw[style=thick] (0.5,2.5) to[out=250,in=135] (0.09,1.4);

\draw[radius=.05, fill=black](0.5,2.5)circle;
\draw[radius=.05, fill=black](1,1)circle;
\draw[radius=.05, fill=black](0,1)circle;

\node at (1.5,1.75) {=};

\end{scope}
%---------10-F2
\begin{scope}[shift={(4,-5.5)}]
\draw[style=thick, dashed](-0.25,1)--(1.25,1);
\draw[style=thick, dashed](-0.25,2.5)--(1.25,2.5);
\draw[style=thick](0.5,2.5)--(0.5,1.75);
\draw[style=thick](1,1)--(0.5,1.75);
\draw[style=thick](0,1)--(0.5,1.75);

\draw[radius=.05, fill=black](0.5,2.5)circle;
\draw[radius=.05, fill=black](1,1)circle;
\draw[radius=.05, fill=black](0,1)circle;
\node at (1.5,1.75) {=};
\end{scope}
%---------10-F3
\begin{scope}[shift={(6,-5.5)}]
\draw[style=thick, dashed](-0.25,1)--(1.25,1);
\draw[style=thick, dashed](-0.25,2.5)--(1.25,2.5);

\draw[style=thick] (0.5,1.75) to[out=180,in=90] (0,1);
\draw[style=thick] (0.5,1.75) to[out=0,in=90] (1,1);
\draw[style=thick] (0.5,2.5) to[out=-70,in=45] (0.92,1.4);

\draw[radius=.05, fill=black](0.5,2.5) circle;
\draw[radius=.05, fill=black](1,1)circle;
\draw[radius=.05, fill=black](0,1)circle;
\end{scope}
\end{tikzpicture}
\caption{Cup and cap diagrams and some isotopy relations in Snyder's diagrammatics.}
\label{fig_1_6}
\end{figure}

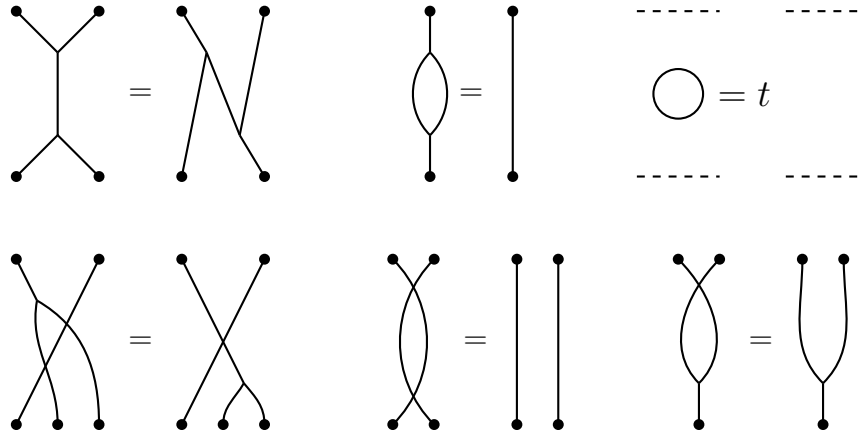
\begin{figure}[h]
\centering
\begin{tikzpicture} [baseline=0cm, scale=1.1]
\draw[style=thick](0,0.5)--(0.5,1);
\draw[style=thick](1,0.5)--(0.5,1);
\draw[style=thick](0.5,2)--(0.5,1);
\draw[style=thick](0.5,1)--(0.5,1);
\draw[style=thick](0.5,2)--(0,2.5);
\draw[style=thick](0.5,2)--(1,2.5);
\draw[radius=.05, fill=black](0,0.5)circle;
\draw[radius=.05, fill=black](0,2.5)circle;
\draw[radius=.05, fill=black](1,0.5)circle;
\draw[radius=.05, fill=black](1,2.5)circle;
\node at (1.5,1.5) {=};
%---------A2
\begin{scope}[shift={(2,0)}]
\draw[style=thick](0,0.5)--(0.3,2);
\draw[style=thick](1,0.5)--(0.7,1);
\draw[style=thick](0.3,2)--(0.7,1);
\draw[style=thick](0.3,2)--(0,2.5);
\draw[style=thick](0.7,1)--(1,2.5);

\draw[radius=.05, fill=black](0,0.5)circle;
\draw[radius=.05, fill=black](0,2.5)circle;
\draw[radius=.05, fill=black](1,0.5)circle;
\draw[radius=.05, fill=black](1,2.5)circle;
\end{scope}

%---------B1
\begin{scope}[shift={(5,0)}]
\draw[style=thick](0,0.5)--(0,1);
\draw[style=thick] (0,1)to[out=135,in=-135] (0,2);
\draw[style=thick] (0,1)to[out=45,in=-45] (0,2);
\draw[style=thick](0,2)--(0,2.5);
\draw[radius=.05,fill=black](0,0.5)circle;
\draw[radius=.05,fill=black](0,2.5)circle;
\node at (0.5,1.5) {=};
\draw[style=thick](1,0.5)--(1,2.5);
\draw[radius=.05,fill=black](1,0.5)circle;
\draw[radius=.05,fill=black](1,2.5)circle;
\end{scope}

%---------C1
\begin{scope}[shift={(8,0)}]
\draw[style=thick, dashed](-0.5,0.5)--(0.5,0.5);
\draw[style=thick, dashed](-0.5,2.5)--(0.5,2.5);
\draw[style=thick,radius=.3](0,1.5)circle;
% \node[cell,label=below:{\Large \i}] at (0,0);
\end{scope}
%---------C2
\begin{scope}[shift={(9.5,0)}]
\node at (-0.7,1.5) {\Large{$=t$}};
\draw[style=thick, dashed](-0.2,0.5)--(0.7,0.5);
\draw[style=thick, dashed](-0.2,2.5)--(0.7,2.5);
\end{scope}

%new line...
%---------D1
\begin{scope}[shift={(0,-3)}]
\draw[style=thick](0,0.5)--(1,2.5);
\draw[style=thick](0,2.5)--(0.25,2);

\draw[style=thick] (0.25,2)to[out=260,in=90] (0.5,0.5);
\draw[style=thick] (0.25,2)to[out=-30, in=90] (1,0.5);

\draw[radius=.05, fill=black](0,0.5)circle;
\draw[radius=.05, fill=black](0.5,0.5)circle;
\draw[radius=.05, fill=black](1,0.5)circle;
\draw[radius=.05, fill=black](0,2.5)circle;
\draw[radius=.05, fill=black](1,2.5)circle;
\node at (1.5,1.5) {=};
\end{scope}
%---------D2
\begin{scope}[shift={(2,-3)}]
\draw[style=thick](0,0.5)--(1,2.5);
\draw[style=thick](0,2.5)--(0.75,1);

\draw[style=thick] (0.75,1)to[out=240,in=90] (0.5,0.5);
\draw[style=thick] (0.75,1)to[out=-50, in=90] (1,0.5);

\draw[radius=.05, fill=black](0,0.5)circle;
\draw[radius=.05, fill=black](0.5,0.5)circle;
\draw[radius=.05, fill=black](1,0.5)circle;
\draw[radius=.05, fill=black](0,2.5)circle;
\draw[radius=.05, fill=black](1,2.5)circle;
\end{scope}

%---------E1
\begin{scope}[shift={(4.55,-3)}]
\draw[style=thick] (0,0.5)to[out=45,in=-45] (0,2.5);
\draw[style=thick] (0.5,0.5)to[out=135,in=225] (0.5,2.5);

\node at (1,1.5) {=};
\draw[style=thick](2,0.5)--(2,2.5);
\draw[style=thick](1.5,0.5)--(1.5,2.5);
\draw[radius=.05,fill=black](0,0.5)circle;
\draw[radius=.05,fill=black](0.5,2.5)circle;
\draw[radius=.05,fill=black](0.5,0.5)circle;
\draw[radius=.05,fill=black](0,2.5)circle;

\draw[radius=.05,fill=black](1.5,0.5)circle;
\draw[radius=.05,fill=black](2,0.5)circle;
\draw[radius=.05,fill=black](2,2.5)circle;
\draw[radius=.05,fill=black](1.5,2.5)circle;
\end{scope}

%---------F1
\begin{scope}[shift={(8,-3)}]
\draw[style=thick] (0.25,1)to [out=45,in=-45] (0,2.5);
\draw[style=thick] (0.25,1)to [out=135,in=225]  (0.5,2.5);
\draw[style=thick](0.25,0.5)--(0.25,1);
\node at (1,1.5) {=};
\draw[style=thick] (1.75,1)to [out=135,in=270](1.5,2.5);
\draw[style=thick] (1.75,1)to [out=45,in=270] (2,2.5);
\draw[style=thick](1.75,0.5)--(1.75,1);
\draw[radius=.05,fill=black](0.25,0.5)circle;
\draw[radius=.05,fill=black](0.5,2.5)circle;
\draw[radius=.05,fill=black](0,2.5)circle;

\draw[radius=.05,fill=black](1.75,0.5)circle;
\draw[radius=.05,fill=black](2,2.5)circle;
\draw[radius=.05,fill=black](1.5,2.5)circle;
\end{scope}
 \end{tikzpicture}
\caption{Some other defining  relations in the Snyder  calculus.}
\label{fig_1_7}
\end{figure}

It is convenient to introduce cup and cap diagrams (as additional generators), defined in the top row of Figure~\ref{fig_1_6} via the original generators. 
Isotopy relations on the generators are shown in the next two rows of  Figure~\ref{fig_1_6}. Some other defining relations are shown in Figure~\ref{fig_1_7}. We leave  it to the reader to convert a full set of relations as found in~\cite[Theorem 1]{C} or~\cite{LS} into defining relations for the Snyder calculus. 

\begin{figure}[!htb]
\centering
\begin{tikzpicture} 
\draw[style=thick](0,0.5)--(0.5,1);
\draw[style=thick](1,0.5)--(0.5,1);
\draw[style=thick](0.5,2)--(0.5,1);
\draw[style=thick](0.5,1)--(0.5,1);
\draw[style=thick](0.5,2)--(0,2.5);
\draw[style=thick](0.5,2)--(1,2.5);
\draw[radius=.05, fill=black](0,0.5)circle;
\draw[radius=.05, fill=black](0,2.5)circle;
\draw[radius=.05, fill=black](1,0.5)circle;
\draw[radius=.05, fill=black](1,2.5)circle;
\node at (1.5,1.5) {=};
%---------A2
\begin{scope}[shift={(2,0)}]
\draw[style=thick](0,0.5)--(0.3,2);
\draw[style=thick](1,0.5)--(0.7,1);
\draw[style=thick](0.3,2)--(0.7,1);
\draw[style=thick](0.3,2)--(0,2.5);
\draw[style=thick](0.7,1)--(1,2.5);

\draw[radius=.05, fill=black](0,0.5)circle;
\draw[radius=.05, fill=black](0,2.5)circle;
\draw[radius=.05, fill=black](1,0.5)circle;
\draw[radius=.05, fill=black](1,2.5)circle;

\node at (2,1.5) {$\Leftrightarrow$};
\end{scope}

\begin{scope}[shift={(5,0.5)}]

 \draw[style=thick] (0.3,2) to[out=270, looseness =2,in=270] (0.7,2);
 
  \draw[style=thick] (-0.3,2) to[out=300, looseness =2,in=70] (-0.3,0);

  \draw[style=thick] (1.3,2) to[out=240, looseness =2,in=110] (1.3,0);
  
\draw[style=thick] (0.7,0) to[out=90, looseness =3,in=90] (0.3,0);

 \draw[fill=white] (0,2) ellipse (0.3cm and 0.1cm);
 \draw[fill=white](1,2) ellipse (0.3cm and 0.1cm);
 
 \draw[fill=white] (0,0) ellipse (0.3cm and 0.1cm);
 \draw[fill=white](1,0) ellipse (0.3cm and 0.1cm);
 \node at (1.5,1) {$\simeq$};
\end{scope} 

\begin{scope}[shift={(7,0.5)}]

%curve top
 \draw[style=thick] (0.3,2) to[out=-120, looseness =1,in=120] (0.4,1.2);
 \draw[style=thick] (0.4,1.2) to[out=-60, looseness =0.8,in=160] (0.7,1);
 
  \draw[style=thick] (0.7,1) to[out=20, looseness =0.8,in=280] (0.7,2);

%curvebottom  
  \draw[style=thick] (0.7,0) to[out=60, looseness =1.1,in=310] (0.6,0.8);
  \draw[style=thick] (0.6,0.8) to[out=110, looseness =0.8,in=340] (0.3,1);  
\draw[style=thick] (0.3,1) to[out=230, looseness =0.8,in=110] (0.3,0);

 \draw[style=thick] (-0.3,2) to[out=300, looseness =1,in=70] (-0.3,0);

  \draw[style=thick] (1.3,2) to[out=240, looseness =1,in=110] (1.3,0);
  
 \draw[fill=white] (0,2) ellipse (0.3cm and 0.1cm);
 \draw[fill=white](1,2) ellipse (0.3cm and 0.1cm);
 
 \draw[fill=white] (0,0) ellipse (0.3cm and 0.1cm);
 \draw[fill=white](1,0) ellipse (0.3cm and 0.1cm);
\end{scope} 
\end{tikzpicture}
\caption{One of the defining relations versus surface diffeomorphism.}
\label{fig_1_8}
\end{figure}
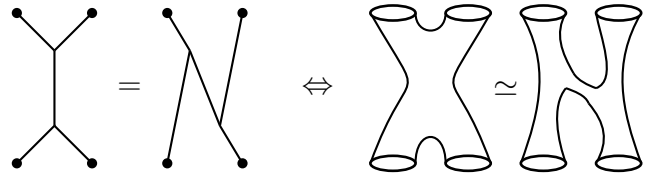

The relation between 2D cobordisms and partitions is especially easy to see in the Snyder calculus. Thickening Snyder's trivalent graphs when  viewed as graphs in $\R^3$ rather than in $\R^2$ results in a surface with boundary that corresponds to the partition. Intersection points of different components  of the graph should be disregarded, as before, for instance by pulling the components slightly  apart in $\R^3$ before thickening (the embedding into $\R^3$ is then  forgotten). An example of a  matching between Snyder's relations and diffeomorphisms of surfaces is shown in Figure~\ref{fig_1_8}. 

\vspace{0.1in}

{\bf From graphs to surfaces.}
The lifting, discussed in this paper and in~\cite{C}, from partitions,  which  are graph-like objects, to  two-dimensional cobordisms is analogous, in some rather naive way, to passing from Feynman diagrams (graph-like objects) to strings (two-dimensional objects): 
\begin{center}
    \begin{tabular}{ccl} 
    &  &  \\
      Feynman diagrams & $\lra$  &   Strings \\ 
        &  &      \\
        Partition diagrams &  $\lra$ & 2D cobordisms \\
        &  & 
    \end{tabular} 
\end{center}
Of course, the complexity of mathematics hidden in the top arrow structures is orders of magnitude higher than those in the bottom arrow, discussed in the present paper. 

\vspace{0.1in}

{\bf Deligne category.}
Let us specialize to ground ring $R=\kk$ a field of characteristic $0$. 
The Deligne category $\delcat$ is the additive Karoubi envelope of the partition category $\Pa_t$, $t\in \kk$, 
\begin{equation}
    \delcat = \mathrm{Kar}(\Pa_t^{\oplus}). 
\end{equation}
It is known to be semisimple when $t\notin \Z_+$. When $t=n\in \Z_+\subset \kk$, the Deligne category admits a nontrivial ideal $J_n$ that consists of \emph{negligible morphisms}. A morphism $x\in \Hom(a,b)$ is negligible if for any $y\in \Hom(b,a)$ the trace of the composition $\Tr(yx)=0$. 
The category $\delcat$ is a tensor category with duals, and the trace is straightforward to define. The trace in $\Pa_t$ on a diagram $\lambda\in D^m_m$ is given by identifying points $i$ and $i'$, $1\le i\le m$. If $r$ is  the number of components in the resulting diagram, $\Tr(\lambda)=t^r$. 

The quotient $\delcatn/J_n$ is equivalent, as a tensor category, to the category $\kk[S_n]\dmod$ of  finite-dimensional representations of the symmetric group $S_n$, 
\begin{equation}
    \delcatn/J_n  \cong \kk[S_n]\dmod . 
\end{equation}
The ideal $J_t$ of negligible morphisms in $\delcat$ is  zero if $t\notin \Z_+$.

%%%%%%%%%%%%%%%%%%%%%%
%
%  Generalized Deligne  
%
%%%%%%%%%%%%%%%%%%%%%%

\section{Generalized Deligne categories} 
\label{sec_gen_Deligne} 

\smallskip 

{\bf Cobordisms  and partitions.}
Consider the  category $\Cobtwo$  of two-dimensional cobordisms.  Given a morphism  $x$ from $n$ to $m$,  disregard its closed components and ignore genera of connected components with boundary. A connected component with boundary defines a subset among the set of boundary circles of $x$. The latter set can be identified with $\{1,2,\dots, n, 1',2',\dots, m'\}$, see Figure~\ref{fig_1_1}.
Consequently, the union of connected components of $x$ that have a  non-empty boundary determines a partition in $D^m_n$. To a cobordism $x$ from $n$ to $m$ we associate this partition in $D^m_n$, denoted $p(x)$. 

To extend this assignment to a functor 
\begin{equation}
    F\ : \ R\Cobtwo \lra  \Pa_t
\end{equation}
let $|cl(x)|$ be the number of connected components of $x$ without boundary (closed components).  Functor $F$ is identity on objects $n\in \Z_+$  of $R\Cobtwo$ and  $\Pa_t$ and 
\begin{equation}\label{eq_func_F}
    F(x) = t^{|cl(x)|} p(x) 
\end{equation}
on cobordisms. It is then extended $R$-linearly to linear combinations of cobordisms. Notice that $F$ ignores genera of all components of $x$. Clearly, $F$ is a tensor (symmetric monoidal) functor. This construction can be found in Comes~\cite[Section 2.2]{C}. One  can   think of $\Pa_t$ as the quotient of  $R\Cobtwo$ by  skein relations in Figure~\ref{fig_2_0}. 

\vspace{0.07in} 

Recall categories $\Cobalp$ and  $\Cobal$ introduced earlier and  associated to a sequence $\alpha$, where a closed surface of genus $g$ evaluates to $\alpha_g\in R$. Let $\alpha(t)=(t,t,t,\dots)$ be the constant sequence associated to $t\in \kk$. The evaluation $\alpha(t)$ associates $t$ to any oriented connected closed surface irrespectively of its genus. Relations in Figure~\ref{fig_2_0} hold in the category $\Cob_{\alpha(t)}$
and they hold in $\Cob'_{\alpha(t)}$ when restricted  to closed components. Consequently, there  are  natural  
tensor functors 
\begin{equation}
 \Cob'_{\alpha(t)} \stackrel{F'_t}{\lra} \Pa_t \stackrel{F''_t}{\lra} \Cob_{\alpha(t)}
\end{equation}
between these three categories. These  functors are 
identities on objects,  $F'_t(n)=n$, $F''_t(n)=n$. The first functor forgets about handles of each component of  a cobordism $S$, evaluates each closed  component  to  $t$, and  associates a partition of the set $\N_n^m$ in (\ref{eq_set_n_m}) to $S$ according to subsets of boundary  circles bounded by connected components of $S$. 

The second functor $F''_t$ exists by an earlier discussion, due to the definition of $\Cob_{\alpha(t)}$ via the quotient by the kernel of a bilinear form. It identifies $\Cob_{\alpha(t)}$ with the quotient of $\Pa_t$ by the ideal of  negligible  morphisms.

\vspace{0.1in} 

{\bf The Deligne category.}
Starting with the functor $F''_t$ and  passing to additive Karoubi closures results in a functor 
\begin{equation} 
F_t \ : \ \delcat \lra \Kob_{\alpha(t)}
\end{equation} 
from the Deligne category to the additive Karoubi closure $\Kob_{\alpha(t)}=\Kar(\Cob^{\oplus}_{\alpha(t)})$ of the category $\Cob_{\alpha(t)}$. From the structure theory of the Deligne  categories we  can conclude that  $F_t$ consists of taking the quotient of $\delcat$ by  the ideal $J_t$ of  negligible morphisms and induces an equivalence
\begin{equation} 
 \delcat/J_t \stackrel{\cong}{\lra} \Kob_{\alpha(t)}.
\end{equation} 
Notice that there's a difference in the order in which we take the additive Karoubi closure and mod out by negligible morphisms. On the Deligne  category side, one  first forms the additive Karoubi closure and  then mods out by negligible morphisms.  On  the $\Kob_{\alpha(t)}$ side, one  first  mods out by negligible morphisms to get the category $\Cob_{\alpha(t)}$ and then  forms the additive Karoubi closure. It is not clear whether this may produce a discrepancy in more general cases, but for Deligne  categories (and with $R$ a field $\kk$ of characteristic $0$) this change of order results in equivalent categories and makes no difference. 

\vspace{0.1in} 

{\bf Generalizations.} We obtain an immediate generalization of 
the categories $\delcat/J_t$ by changing from the constant sequence $\alpha(t)$ in (\ref{seq_alpha_t}) to a more general sequence $\alpha$. The most interesting case is when the generating function $Z_{\alpha}(T)$ of $\alpha$, see (\ref{eq_Z_pow}), is a rational function, a ratio of two coprime polynomials
\begin{equation}\label{eq_Z_alpha}
    Z_{\alpha}(T) = \frac{P(T)}{Q(T)}
\end{equation}
with coefficients  in $\kk$. In  this case categories $\Cobal$ and $\Kobal$ have finite-dimensional hom spaces. We can view $\Kobal$ as a natural generalization of the quotient category $\delcat/J_t$. For generic $t$, the ideal $J_t$ is zero, and then  the quotient category is the Deligne category. 

\begin{theorem}
Categories $\Kobal$ are tensor $\kk$-linear Karoubi-closed additive categories. When $Z_{\alpha}(T)$ is rational,    morphism spaces in $\Kobal$  are 
finite dimensional. 
\end{theorem} 

It is an interesting  project to investigate categories $\Kobal$ when the generating function $Z_{\alpha}(T)$ is rational. Deligne category quotients are recovered for the rational function in  (\ref{eq_D_rat_fun}). 

\vspace{0.1in} 

Notice that categories $\Kobal$ deliver  generalizations of the quotients 
 $\delcat/J_t$ rather than of Deligne categories $\delcat$ themselves. To remedy this discrepancy, we instead pass from $\Cobalp $ to $\Kobal$ in one more step, when $R=\kk$ is a  field and the  partition  function $Z_{\alpha}(t)$ is rational (\ref{eq_Z_alpha}).  Let 
 \begin{eqnarray} \label{eq_K}
  N & = & \deg(P(T)), \ \ M=\deg(Q(T)), \ \  K = \max(N+1,M), \\
  Q(T) & = & 1 - e_1 T +e_2 T^2  + \ldots + (-1)^M e_M T^M, \ e_i \in \kk, 
 \end{eqnarray}
 as in~\cite[Section 2.4]{Kh2}. Then in the state space $A_{\alpha}(1)$ of a circle equality  
 \begin{equation}\label{eq_with_K} 
     x^K - e_1 \, x^{K-1}+e_2 \, x^{K-2}- \ldots + (-1)^M e_M \, x^{K-M} =0 
 \end{equation}
 holds, where  $x$ denotes a 2-torus with one boundary component. Power $x^k$ of $x$ represents  a surface of genus $k$ with one boundary component, with multiplication in $A_{\alpha}(1)$ given by the pants cobordism, see~\cite{Kh2}.   
Equation (\ref{eq_with_K}) gives a skein relation in category $\Cobal$ which reduces a collection of $K$ handles on a single component to a linear combinations of collections of $K-1$, $K-2$, ... , $ K-M$ handles. 

For rational $\alpha$, start with the pre-additive category $\Cobalp$, see Section~\ref{sec_cat_lin}, and diagram (\ref{eq_many_cats}) that shows the position of $\Cobalp$ in the chain of categories and functors associated with $\alpha$. In $\Cobalp$ only closed components are reduced to elements $\alpha_k$ of  $\kk$. Hom spaces $\Hom(n,m)$ in  $\Cobalp$ are infinite-dimensional $\kk$-vector spaces, unless $n=m=0$, with a basis of diffeomorphism classes rel boundary of all cobordisms without closed components. Thus, a basis element is described by a decomposition in $D^m_n$ and a choice of genus for each connected component. 
   
\vspace{0.1in}   

Define category  $\PCobal$ to have the same objects $n\ge 0$ as $\Cobalp$ and morphism spaces to be quotients of those in $\Cobalp$ by the skein relations corresponding to the equation (\ref{eq_with_K}). That is, we set this linear combination of morphisms (cobordisms) to zero in  the quotient category. Applying this relation we reduce a component which contains at least $K$ handles to components with fewer handles. In particular, any  morphism in $\Cobalp$ reduces to a  $\kk$-linear  combination of cobordisms with  no  closed components  and  at most $K-1$ handles on each connected component. Diffeomorphism classes rel  boundary of  these cobordisms  
are  in  a bijection with elements of  the set  $D^m_n(<K)$ of partitions in $D^m_n$ with an integral weight  between $0$ and $K-1$ associated to each part (number of handles of the component, on the cobordism  side). Recall that to a partition $x$ we associated cobordism $p(x)$, see discussion preceeding formula (\ref{eq_func_F}). We can now extend this association, also denoted $p$, and assign to a partition $x$  with non-negative integral weights of its parts the cobordism $p(x)$ by starting with the cobordism for the partition without weights and adding the  number of  handles equal to the weight to each two-sphere with boundary holes. 

Computations in~\cite[Section 2.4]{Kh2} imply that relation (\ref{eq_with_K}) is compatible with evaluation $\alpha$ applied to closed  cobordisms. In  particular,  no additional relations on cobordisms appear and  elements  of the  set  $D^m_n(<K)$, converted to cobordisms, provide a basis of $\Hom_{\PCobal}(n,m)$.

\begin{prop} The hom space $\Hom(n,m)$ in $\PCobal$ has a basis $\{p(x)\}$, over all $x\in D^m_n(<K)$. 
\end{prop} 

In particular, hom spaces in $\PCobal$  are finite-dimensional. 
We can now insert category $\PCobal$ into the chain of six categories in (\ref{eq_many_cats}):
\begin{equation}\label{eq_many_cats_2} 
    \Cobtwo \lra R\Cobtwo\lra \Cobalp \lra \PCobal \lra \Cobal\lra \Cobal^{\oplus} \lra \Kobal.
\end{equation}
It fits in between $\Cobalp$ and $\Cobal$. Category $\PCobal$ is the quotient of $\Cobalp$ by the skein relation (\ref{eq_with_K}). Like every other category  in  this chain, it is tensor (symmetric monoidal). The trace form on $\Cobalp$
descends to that on $\PCobal$. The quotient of $\PCobal$ by the ideal of  neglible morphisms relative to this trace form is isomorphic to $\Cobal$ (isomorphic and not just equivalent, since these categories are essentially skeletal and have very few objects). As we've mentioned, this insertion is possible  when $Z_{\alpha}(T)$ is a  rational function   and  $R$ is a field. 

\vspace{0.1in} 

Category  $\PCobal$ generalizes the partition category $\Pa_t$. Partition category $\Pa_t$ is isomorphic to $\mathrm{PCob}_{\alpha(t)}$ for the constant sequence $\alpha(t)=(t,t,\dots)$.  
More generally, choosing $K$ in (\ref{eq_K}) fixes  the size of  homs in $\PCobal$, analogously to independence of dimensions of homs in $\Pa_t$ on $t$, given by the number of partitions (the  Bell number). When $K=1$, dimensions of hom spaces in $\PCobal$ are also given by the number of partitions in $D^m_n$, and for $K>1$ by the number of weighted partitions as discussed above. 

\vspace{0.1in} 

To get from  $\PCobal$ to the analogue of the Deligne category, pass to the additive Karoubi closure to get an additive and idempotent-complete  category 
\begin{equation}
    \KPobal \ := \ \Kar(\PCobal^{\oplus}). 
\end{equation}

Chain of  functors (\ref{eq_many_cats_2}) can be  upgraded to a commutative diagram of functors

\begin{equation} \label{eq_seq_cd}
\begin{CD}
\Cobtwo @>>> R\Cobtwo @>>> \Cobalp @>>> \PCobal  @>>> \PCobal^{\oplus} @>>> \KPobal  \\
@.  @.  @.   @VVV      @VVV  @VVV  \\    
 @.   @.    @.  \Cobal @>>> \Cobal^{\oplus} @>>> \Kobal
\end{CD}
\end{equation}
where the chain (\ref{eq_many_cats_2}) is given by  the left, left, left, down, left, left sequence of arrows. Two new categories are added in the upper right. Vertical down arrows are  quotients by the ideals of negligible morphisms. Both squares in the diagram is commutative. 

Notice that first  modding out $\PCobal$ by  negligible morphisms to get $\Cobal$ and then taking the Karoubi envelope $\Kobal$ compared to first taking the Karoubi envelope  $\KPobal$ and then modding out by  negligible morphisms does not produce any extra idempotents. This is due to the easy to check idempotent lifting property that holds for any finite-dimensional algebra $B$ over $\kk$ and any 2-sided ideal $J\subset B$ (not necessarily nilpotent). Any idempotent in $B/J$ lifts to an idempotent in $B$. Endomorphism algebras of objects in $\PCobal^{\oplus}$ are finite-dimensional over $\kk$. For the  ideal $J$ one would take the  ideal of negligible endomorphisms of an object in  $\PCobal^{\oplus}$.

\vspace{0.1in}

To summarise, the chain of three  categories  and  two functors (the partition category, the Deligne category, and its quotient by  negligible morphisms) 
\begin{equation}\label{eq_three_cat}
    \Pa_t \lra \delcat \lra \delcat/J_t 
\end{equation}
generalizes to a similar chain 
\begin{equation}\label{eq_many_cats_4} 
     \PCobal \lra \KPobal \lra  \Kobal
\end{equation}
for any  sequence $\alpha$ with  rational power series $Z_{\alpha}(T)$. Specializing  to the constant series $\alpha(t)$ and  rational function $t/(1-T)$ recovers the  original setup (\ref{eq_three_cat}). 

\medskip 

%%%%%%%%%%%%%%%%%%%%%
%%
%%   REFERENCES 
%%
%%%%%%%%%%%%%%%%%%%%

%\addcontentsline{toc}{section}{References}
%\def\refname{}

\end{document}